\documentclass[11pt,a4paper]{article}

\usepackage{cite}
\usepackage{amsmath,bm}
\usepackage{amsthm}
\usepackage{amsfonts}
\usepackage{amssymb}
\usepackage{latexsym}
\usepackage[T1]{fontenc}
\usepackage[utf8]{inputenc}
\usepackage{csquotes}
\usepackage{caption}
\usepackage{graphicx,subcaption,fancyhdr,multirow,stmaryrd}
\usepackage[height=22 cm,width=16 cm,top=4 cm,left=3 cm]{geometry}
\usepackage{epstopdf}
\usepackage{wrapfig} 

\usepackage{xcolor}

\graphicspath{{figures/}}

\def\N{\mathbb{N}}

\def\R{\mathbb{R}}

\def\e{\mbox{e}}
\def\d{\mbox{d}}

\def\ul#1{{\mathbf #1}}

\def\s{\textstyle}

\def\e{\text{e}}

\bmdefine\bnu{\nu}
\bmdefine\bmu{\mu}
\bmdefine\bmeta{\eta}

\numberwithin{equation}{section}

\theoremstyle{definition}

\theoremstyle{remark}

\begin{document}


\title {On conservative, stable boundary and coupling conditions for diffusion equations I - The conservation property
for explicit schemes}

\author{
Taj Munir\thanks{\baselineskip=.5\baselineskip Department of Energy and Power Engineering, 
Jiangsu University, Zhenjiang, 212013, China.
Email: taj\underline{ }math@hotmail.com},
Nagaiah Chamakuri\thanks{\baselineskip=.5\baselineskip School of Mathematics,
Indian Institute of Science Education and Research (IISER-TVM), Maruthamala PO, Vithura
Thiruvananthapuram, Kerala, India-695551.
Email:  nagaiah.chamakuri@iisertvm.ac.in}
and
Gerald Warnecke\thanks{\baselineskip=.5\baselineskip Institute of Analysis and  Numerics,
Otto-von-Guericke-University Magdeburg, PSF 4120, D--39016 Magdeburg,
Germany.  Email: gerald.warnecke@ovgu.de}
}

\date {\today }

\maketitle

\begin{abstract}
This paper introduces improved numerical techniques for addressing numerical boundary and interface coupling conditions 
in the context of diffusion equations in cellular biophysics or heat conduction problems in fluid-structure interactions. 
Our primary focus is on two critical 
numerical aspects related to coupling conditions: the preservation of the conservation property and ensuring stability. 
Notably, a key oversight in some existing literature on coupling methods is 
the neglect of upholding the conservation property 
within the overall scheme. This oversight forms the central theme of the initial part of our research. 
As a first step, we limited ourselves to
explicit schemes on uniform grids. Implicit schemes and the consideration of varying mesh sizes at the interface will be reserved for a 
subsequent paper \cite{CMW3}. Another paper \cite{CMW2} will address the issue of stability.

We examine these schemes from the perspective of finite differences, including finite elements, 
following the application of a 
nodal quadrature rule. Additionally, we explore a finite volume-based scheme involving cells 
and flux considerations. Our 
analysis reveals that discrete boundary and flux coupling conditions uphold the conservation property 
in distinct ways in nodal-based and cell-based schemes.
The coupling conditions under investigation encompass well-known approaches such as Dirichlet-Neumann coupling, 
heat flux coupling, 
and specific channel and pumping flux conditions drawn from the field of biophysics. The theoretical findings pertaining to 
the conservation property are corroborated through computations across a range of test cases.  
\end{abstract}



\section{Introduction}
\label{sec:int}

The subject of this paper is a study of various numerical interface coupling conditions for 
diffusion or heat equations.     
Many important real-world problems in physics, engineering, and biology are modeled via bi-domain or multi-domain
partial differential equations (PDEs) with coupling conditions at the sub-domain interfaces. 
This modeling may be due to the nature of the problem when it has a
physical interface at which physical coupling conditions occur. Or it may be an artificial mathematical interface, 
e.g.\ in domain decomposition methods that are used for computational purposes. 


Diffusion equations, commonly used to model heat conduction via Fourier's law or the diffusion of mass concentrations via Fick's law, are central to 
many of these models. One such application is the diffusion of calcium concentrations in living cells, which is modeled using a system of reaction-
diffusion equations. A complex three-dimensional biophysical model, developed by Falcke \cite{l8}, incorporates coupling conditions to describe 
intracellular calcium dynamics between the cytosolic and endoplasmic reticulum (ER) regions of a cell. The channels and pumps on the membrane 
separating these regions are mathematically modeled as coupling conditions on the interface between them, 
as explored by Thul \cite{l10}, Thul and Falcke \cite{l11}, as well as 
Chamakuri \cite{l9}. These coupling conditions were the main motivation for this study.
Another motivation came from the analysis of Giles \cite{GIL} in the context of fluid-structure interactions.


While the diffusion equations need only one boundary condition at outer boundaries,
there are always two conditions needed for coupling them at an internal interface. 
The coupling conditions of the above model are channel pumping and membrane pumping fluxes. 
For the purpose of comparison, we also consider some of the other types of coupling conditions used in conjunction with
bi-domain diffusion or heat equations.
These include the well-known Dirichlet-Neumann coupling and the heat flux coupling. 

The simplest coupling conditions are the Dirichlet-Neumann (DN) coupling conditions that
have been extensively used for parallel computing in non-overlapping domain decomposition methods, see e.g.\ 
Toselli and Widlund \cite{bTOWI}, Quarteroni and Valli \cite{bQUVA} or Quarteroni \cite[Ch.\ 19]{bQUA}.
The same type of problems also arise in heat flow at interfaces between different materials 
where heat flux coupling conditions are well
established, see e.g.\ Carslaw and Jaeger \cite{b8}. 
Carr and March \cite{l116} considered various interface coupling conditions based on the heat flux coupling conditions. 
The coupling conditions in which we are mainly interested are the more complex channel pumping and membrane
pumping conditions modeling calcium transport within cells mentioned above, 
see these conditions in \eqref{ch66} and \eqref{987} below. 
They are extensions of the heat flux coupling conditions, e.g.\ by the addition of a non-linear
pumping term. 
  
Our analysis in this first paper is focused on the essential property of conservation of concentration or heat.
This property reflects the fundamental principles of mass conservation for mass concentration or energy conservation
for temperature.
Indeed, this property must be maintained at the discrete level by numerical schemes. We will show
how this is relevant to discretized flux boundary conditions and various
flux coupling conditions for bi-domain diffusion equations. As a first step, we only look at schemes that
are explicit in time. We are aware of the fact that implicit methods are
actually the methods of choice for these stiff problems. The implicit schemes are addressed in a further paper \cite{CMW3}.
In a companion paper \cite{CMW2}, we will consider the property of 
Godunov-Ryabenkii stability is based on normal mode solutions for the coupling schemes as
well as for the boundary conditions. Results of the joint research of the authors 
were included in the thesis of Munir \cite{MUN}.
Until now, the bio-physical coupling conditions that we are considering
have not been analyzed in terms of the discrete conservation property and numerical stability.

The use of GR stability for coupling conditions was first introduced by Giles \cite{GIL}, which served as a major starting point for our studies. 
However, Giles's approach employed an inconsistent scheme, leading to artificial instabilities and loss of conservation. This inconsistency was 
pointed out by Zhang et al.\ \cite{zhang}. We will address it in our discussion of the schemes. Giles focused on one-dimensional bi-domain 
thermal diffusion equations with Dirichlet-Neumann coupling conditions, which he used in engineering heat conduction problems involving fluid-
structure interactions. For simplicity, we consider only the basic diffusion equation.

Roe et al.\ \cite{roe2007stability} extended Giles's work by introducing a moving interface, using both finite difference/finite difference (FD/FD) and finite volume/finite element (FV/FEM) discretizations. In a later study \cite{roe2008combined}, they used higher-order combined methods for explicit and implicit coupling. Errera and Chemin \cite{errera2013optimal} explored Dirichlet-Robin and Robin-Robin coupling conditions for the same equations, and Errera with Moretti et al.\ \cite{moretti2018stability} focused on stability and convergence for various coupling schemes.

Henshaw and Chand \cite{chand} studied a heat transfer problem as a multi-domain issue with Dirichlet-Neumann coupling and mixed Robin conditions, proposing the use of central differences for discretizing interface equations to improve accuracy and stability. Lemari\'{e} et al. \cite{lemarie2015analysis} explored ocean-atmosphere coupling conditions using implicit and explicit methods, while Zhang et al.\ \cite{zhang} investigated multi-domain PDEs in climate models, employing both explicit and implicit coupling discretizations.

For spatial discretization, we use the standard second order central difference for the second derivative and limit ourselves to a first-order 
method in time. For explicit methods, we employ a simple forward Euler time step, resulting in the forward in time, central in space (FTCS) update. This 
discretization of diffusion equations comes with a significant stability restriction, making it a stiff problem. Consequently, implicit methods are 
generally more practical. However, due to the length of this paper, implicit schemes will be discussed in a separate publication \cite{CMW3}. For basic 
numerical methods and concepts related to single-domain diffusion equations, we refer to Hundsdorfer and Verwer \cite{bHUVE}, Morton and Mayers \cite{bMOMA}, or Thomas \cite{bTHO}.
 
As outer boundary conditions, we are interested in the more commonly used flux conditions rather than the Dirichlet conditions
used by Giles \cite{GIL}. For simplicity, we took
homogeneous Neumann conditions. They have the advantage that total concentration should be conserved on the domain 
and we can test that numerically.
Non-zero fluxes can be easily added. 
Ghost values are used to determine the numerical updates of the outer boundaries and the interface coupling 
conditions. To find these ghost values, we use either the central difference method or a one-sided
difference method with respect to the mesh points at or next to the boundary of the computational domain. 
The ghost values are thereby eliminated from the updates.

Below in Section \ref{sec:cons}, we will argue that it is important to maintain the conservation property in discretizations.
We emphasize that it is necessary to maintain e.g.\ the mass conservation principle for a concentration
in the type of applied problems that we are considering. 
This is the analogue of energy conservation in case the same type of equations model heat flow. 
The numerical conservation property that we introduce restricts the discretization of boundary and coupling conditions.
This is the first important result of the paper. 

The next key result comes from the comparison of nodal based finite element type and cell based finite volume type schemes. 
For the Dirichlet-Neumann coupling there is no difference in the schemes due to
the Dirichlet coupling condition. But, for Neumann type boundary conditions and
flux coupling conditions, the schemes differ. Nodal based schemes need a central
difference with respect to the boundary or interface node
in order to have the conservation property.
In finite volume type schemes, which are cell based, there are no nodal values 
at boundaries or interfaces. The domain boundary or coupling interface is a cell boundary. Here, a one-sided difference is actually a central difference 
for the domain boundary and it is conservative. This distinction is the second main result of the paper.
We also show that the Giles coupling 
\cite{GIL} fails to satisfy the conservation property. 

One of our test cases has discontinuous data at the coupling interface.
In this case, we show that for the Dirichlet-Neumann coupling the nodal based scheme produces a much larger 
computational error in the total mass concentration than the finite volume type scheme. This is due to the fact that the former has a node at the
interface.


The paper is organized as follows. Section \ref{sec:bi} begins with a 
brief review of the biophysical application that motivated this study, providing the necessary background for understanding the bi-domain diffusion problem and 
the associated coupling conditions. We then define the specific bi-domain diffusion problem under investigation and introduce various coupling 
conditions, some drawn from the literature. In Section \ref{sec:exp}, we present our explicit time-
stepping scheme and the nodal discretization of the equations, along with the numerical boundary and coupling conditions. 
We include those suggested by 
Giles \cite{GIL}. We also introduce a finite volume type scheme. Section \ref{sec:cons} is the core of the paper. We introduce the concept of discrete conservation property for both the nodal-based and finite volume type schemes. This is followed by the results discussed earlier. Finally, Section \ref{sec:tests} presents numerical computations that demonstrate the coupling as well as boundary conditions and verify the conservation property in practice.

\section{Bi-domain modeling and coupling conditions}
\label{sec:bi}

The main motivation for this paper is the interest in coupling conditions arising 
in the mathematical modeling of calcium dynamics in living cells.
In a living cell, calcium $Ca^{2+}$ is transported through channels by pumps, and it diffuses into the cytosol
as well as into the endoplasmic reticulum (ER), and it reacts with buffers. The ER and the cytosol are
separated by a membrane that contains channels through which calcium is exchanged.
The calcium concentration in the ER is denoted by $E$, and in the cytosol by $c$. 
A three-dimensional calcium dynamics model for these processes was proposed by Falcke \cite{l8}, see also 
Chamakuri \cite{l9}, Thul and Falcke \cite{l11} as well as Thul \cite{l10} for more detailed descriptions 
of the mathematical model.

\begin{wrapfigure}{l}{0pt}
\includegraphics[width=3 in]{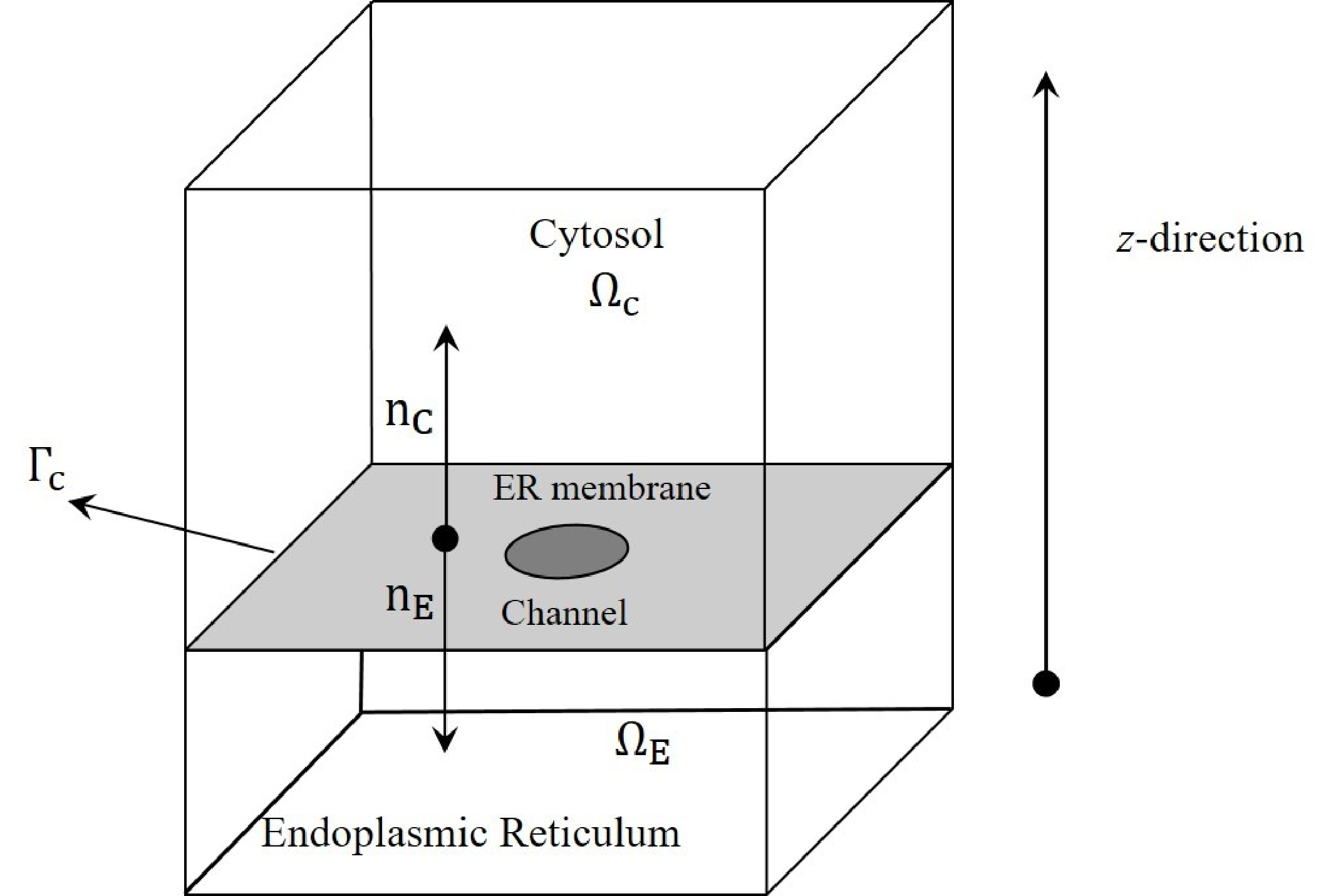}
\centering
\caption{Bi-domain cubic volume distribution of ER and cytosolic domains, modification of a figure from \cite{l9}.}
\label{fig901}
\end{wrapfigure}

In the model, the membrane is a surface that divides
the model domain into two subdomains $\Omega_c$ and $\Omega_E$. When modeling a small section of a cell,
the simplest domain for such a model
is to take a rectangular box or cuboid with a planar surface for the membrane, see Figure~\ref{fig901}.
The calcium dynamics is described
by a system of coupled reaction-diffusion equations for the concentrations of calcium $Ca^{2+}$ in the cytosol, the ER,
 and a number of buffers, see Falcke \cite{l8}. The equations are coupled through reaction terms, which we will disregard. Instead, we will focus on the bi-domain coupling via the membrane. To simplify, we will consider a one-dimensional approximation of the three-dimensional model, as shown in Figure~\ref{fig901}, by taking a one-dimensional slice along the $z$-axis.

\subsection{Biophysical coupling conditions}
\label{subsec:biocouple}
The bi-domain coupling occurs through fluxes that represent channel flow or pumping. We describe these in some detail as we aim to incorporate similar terms into our one-dimensional model. The transport across the ER membrane involves three distinct contributions. Calcium is moved
from the ER into the cytosol through a leak current $P_l (E-c)$ in a pumping term and some channels.
Here, $P_l$ is the coefficient of the leak flux density. Calcium is re-sequestered into the ER by 
pumps modeled by a term proportional to the maximal pump strength $P_p$. The term contains a dissociation constant $K_d>0$. 
In the current model, the channel cluster fluxes $J_{ch}$ and the membrane pumping fluxes $J_{pump}$ are 
normal to the surface, as shown in Figure \ref{fig901}. They are given as follows  \cite[(2.1)]{l10}
\begin{equation}
\label{ch66}
J_{ch}(c,E)=\Psi \frac{E-\alpha c}{\beta+\gamma E+\delta c}
\end{equation}
with some constants $\alpha$, $\beta$, $\gamma$, $\delta$ and $\Psi$. Outside of the channels on the membrane, one has
\begin{equation}
\label{987}
J_{pump}(c,E)=P_{l}(E-c) -P_{p}\frac{c^{2}}{K_{d} ^2 +c^2}
\end{equation}
as membrane pumping condition with some constants $P_l$, $P_p$ and $K_d$.
Specific values of these parameters can be found in Thul \cite{l10} or Thul and Falcke \cite{l11}.

The model consists of reaction-diffusion type equations for $c$ and $E$ as well as a number of other quantities, 
see e.g.\ Thul \cite[Chapter 2]{l10} or Thul and Falcke \cite{l11}.
The respective diffusion coefficients are $D$ and $D_E$.
The currents are incorporated into the volume dynamics by setting the 
flux type coupling conditions on the interface $\Gamma_c$ at the $ER$
membrane to be
\begin{equation}
\label{couple_bio}
D\nabla c \cdot \mathbf{n}_c= D_E \nabla E \cdot \mathbf{n}_c=-J_{ch,pump}(c,E) 
\end{equation}
or more specifically $D\frac{\partial c}{\partial z}=D_E \frac{\partial E}{\partial z}=-J_{ch,pump}(c,E)$.

In Fourier's or Fick's law, the flux $J$ is always in the direction of the negative gradient of the diffused quantity. 
Therefore, a positive slope of the solution must correspond to a negative flux. In the coupling conditions, this means that the 
gradients or normal derivatives of the solution are equal to minus the fluxes.
Note that, in practice, the dynamics of this model are further complicated by the spontaneous and stochastic opening and closing of channels, see Falcke \cite{l8}. However, this is not pertinent to our objectives here.

\subsection[Bi-domain diffusion equation]{The one dimensional bi-domain diffusion equation}
%
%
 
We now consider a bi-domain one-dimensional diffusion model, which is a one-dimensional reduction of the three-dimensional model discussed above. We assume that the concentration only varies in the vertical $z$-direction
and use $x$ as our one-dimensional variable. If it suffices to restrict ourselves to one diffusion equation in order to study 
the coupling conditions. We consider as domain the interval $\Omega=[0,1]\subset\R$ and
divide the domain into two sub-domains by the midpoint $x=\frac 12$, which is the common
interface boundary. The two sub-domains are $\Omega_-=[0,\frac 12]$  and $\Omega_+=[\frac 12, 1]$. 
Let $D_\pm>0$ be the diffusion coefficients on the sub-domains, which may differ. 
For the consideration of the coupling conditions, we want to distinguish the solutions clearly
on the sub-domains. 
We therefore take
$u:\Omega_-\times\mathbb{R}_{\ge 0}\to \mathbb{R}$ and $v:\Omega_+\times\mathbb{R}_{\ge 0}\to \mathbb{R}$ 
to be the solutions that describe a concentration or temperature at
position $x\in\Omega$ and time $t\in\R_{\ge 0}$ in the sub-domains.
We provide some initial data $u_0:\Omega_-\to \mathbb{R}$
and $v_0:\Omega_+\to \mathbb{R}$. The initial boundary value problem for the bi-domain diffusion equations is then defined as
\begin{flalign}
\label{eqn1}
&\quad\frac{\partial u}{\partial t}= D_{-}\frac {\partial ^2 u}{\partial x^2} 
\quad \mbox{for all} \quad (x,t) \in \Omega_+\times\mathbb{R}_{\ge 0},
\qquad\frac{\partial v}{\partial t}= D_{+}\frac {\partial ^2 v}{\partial x^2} 
\quad \mbox{for all} \quad (x,t) \in \Omega_-\times\mathbb{R}_{\ge 0},\nonumber\\
&\quad u(x,0) = u_0(x)\quad\text{for}\;x\in\Omega_-,\qquad\qquad\qquad\qquad v(x,0) = v_0(x)\quad\text{for}\; x\in\Omega_+,\nonumber\\
&\quad\frac{\partial}{\partial x}u(0,t)=\frac{\partial}{\partial x}v(1,t) = 0 \quad\text{for}\; t\in\R_{\ge 0}.
\end{flalign}
We take the outer boundary conditions to be the homogeneous no flux
Neumann conditions. We aim to solve a well-posed problem by coupling $u$ and $v$ across the interface at $x=\frac 12$. The problem
defined in \eqref{eqn1} also needs two appropriate internal coupling conditions at this interface. 

\subsection{Coupling conditions}
\label{subsec:couple}

We now introduce a number of useful coupling conditions that can be found in the literature.

\noindent
{\bf Dirichlet-Neumann coupling:}
The simplest case is to assume the continuity of the solution and the flux at the interface.
\begin{equation}
\label{nnm}
u({\s\frac 12},t)=v({\s\frac 12},t),\qquad D_- 
\frac{\partial u(\frac 12,t)}{\partial x} =D_+ \frac{\partial v(\frac 12,t)}{\partial x}.
\end{equation}
It is generally called the Dirichlet-Neumann coupling.
This type of coupling is used in domain decomposition methods, see e.g.\ Quarteroni and Valli in \cite{bQUVA}. 
In heat conduction, it applies to different materials in perfect contact, as described by Carslaw and Jaeger \cite[p.\ 23]{b8}. Carr and March \cite{l116} refer to it as the perfect contact condition.

Note that in case $D_-=D_+=D$, this coupling will give a solution $w$ to
the single domain diffusion equation. The
common factor $D$ then drops out of \eqref{nnm}. This can be used as a numerical test case for
coupling algorithms.

\noindent
{\bf Heat flux coupling conditions:}
We assume that the heat flow is proportional to the temperature difference and flowing from higher to lower temperature.
\begin{equation}
D_- \frac{\partial u(\frac 12,t)}{\partial x}=D_+ \frac{\partial v(\frac 12,t)}{\partial x} 
=-J_{heat}(u,v) =H(v({\s\frac 12},t)-u({\s\frac 12},t)),
\label{eq111}
\end{equation}
see e.g.\ Carslaw and Jaeger \cite[p.~23]{b8}. Here $H>0$ is the contact transfer coefficient at $x=\frac 12$. 

\noindent
{\bf General interface conditions:}
For completeness, we mention a more general case than the previous coupling conditions.  
All of the above coupling conditions can be considered within a general form as 
\begin{equation}
\label{flux_gen}
D_- \frac{\partial u(\frac 12,t)}{\partial x} =D_+ \frac{\partial v(\frac 12,t)}{\partial x} =-J_{gen}(u,v)
\end{equation}
with $J(u,v)=J_{gen}(u,v) =-H(\theta v({\s\frac 12},t)-u({\s\frac 12},t))$. Here $\theta>0$ is the
partition coefficient at $x=\frac 12$. For $\theta =1$, we have the heat flux coupling conditions \eqref{eq111}.

Using one of the flux equations, dividing it by $H$ and taking $H \to \infty$ produces 
the {\bf partition conditions} given by Carr and March \cite{l116} 
\begin{equation}
u({\s\frac 12},t)=\theta v({\s\frac 12},t),\qquad D_- \frac{\partial u(\frac 12,t)}{\partial x} 
=D_+ \frac{\partial  v(\frac 12,t)}{\partial x}.
\label{part1}
\end{equation}
for $t>0$. The condition defined in case $\theta\neq 1$ maintains a constant
ratio between the discontinuous solutions at the interface. For references to applications, see Carr and March \cite{l116}.
If $\theta =1$, the coupling conditions are just the Dirichlet-Neumann coupling \eqref{nnm}. They are the limiting case for $H\to\infty$ of the heat flux coupling conditions.

\noindent
{\bf Channel and pumping interface conditions:}
Now, we want to consider the coupling conditions \eqref{couple_bio} using
\eqref{ch66} or \eqref{987}.
The general combined coupled interface conditions for the biophysical model discussed above are defined 
using $J(u,v)=J_{ch,pump}(u,v)$ in \eqref{flux_gen} for channel and membrane pumping. The latter flux is given by \eqref{ch66} or \eqref{987} for $u=E$ and $v=c$. 
Note that both cases can be viewed as generalizations of the heat flux conditions \eqref{eq111}.



Note that most coupling conditions have the same form \eqref{couple_bio} or \eqref{flux_gen}. The three terms have to be equal. 
For practical use, one chooses two out of the three possibilities of equating a pair of these terms. These are then discretized and used
in the numerical scheme.
\section{Explicit discretization of the bi-domain diffusion problems}
\label{sec:exp}

In this section, we discuss numerical methods for the discretization of the bi-domain diffusion 
equation with homogeneous Neumann boundary conditions in one space dimension. We discretize the diffusion equation with 
an explicit nodal-based finite difference scheme and a finite volume discretization method.

\subsection{The explicit nodal based scheme}
\label{sub:node}

We take $\Omega=[0,1]$ partitioned into $\Omega_- =[0,\frac 12]$ and $\Omega_+ =[\frac 12, 1]$
with the coupling boundary at $x=\frac{1}{2}$. 
We introduce grid points $x_j\in [0,1]$ for the spatio-temporal discretization of the bi-domain diffusion system. 
We set $N=2m$ for some $m\in\mathbb{N}$ and $\Delta x=1/N$. Grid points for the two sub-domains $\Omega_-=[0,1/2]$
and $\Omega_+=[1/2,1]$ are defined as
\begin{equation*}
x_{j}=j \Delta x,\quad\mbox{for}\quad j=0,1,...,m-1, j=m+1,..., N=2m.
\end{equation*}
The nodes $x_0$ and $x_N$ are the boundary nodes, the others the interior nodes.
Further, at the interface $c=\frac{1}{2}$ we introduce at $x_m =m\Delta x$ a double node $x_{m_{-}}=x_{m_{+}}=m\Delta x$, see
Figure \ref{bid1} below. The node $x_{m_{-}}$ is used 
in conjunction with $\Omega_-$ and $x_{m_{+}}$ with $\Omega_+$. We call a scheme based on this mesh a nodal-based scheme, in contrast to the cell-based finite volume mesh introduced below.
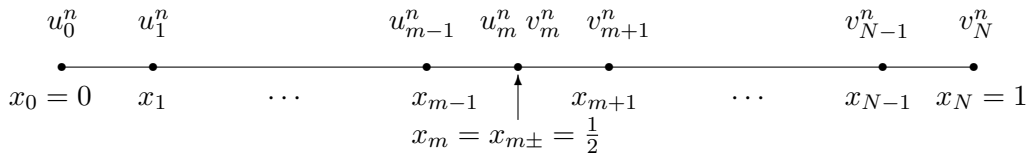
\begin{figure}[htp]
\setlength{\unitlength}{1cm}
\begin{picture}(15,2) 
\put(2.0,1){\line(1,0){12}}       
\put(2.0,1){\circle*{0.1}}
\put(1.3,0.5){$x_0=0$}
\put(1.8,1.5){$u_0^n$}
\put(3.2,1){\circle*{0.1}}
\put(3,0.5){$x_1$}
\put(3,1.5){$u_1^n$}
\put(4.7,0.5){$\cdots$}
\put(6.8,1){\circle*{0.1}}
\put(6.6,0.5){$x_{m-1}$}
\put(6.3,1.5){$u_{m-1}^n$}
\put(8,1){\circle*{0.1}}
\put(6.6,0){$x_m =x_{m\pm}=\frac 12$}
\put(8,0.3){\vector(0,1){0.6}}
\put(7.5,1.5){$u_m^n$}
\put(8.1,1.5){$v_m^n$}
\put(9.2,1){\circle*{0.1}}
\put(8.7,0.5){$x_{m+1}$}
\put(8.9,1.5){$v_{m+1}^n$}
\put(10.8,0.5){$\cdots$}
\put(12.8,1){\circle*{0.1}}
\put(12.3,0.5){$x_{N-1}$}
\put(12.3,1.5){$v_{N-1}^n$}
\put(14,1){\circle*{0.1}}
\put(13.5,0.5){$x_N=1$}
\put(13.8,1.5){$v_N^n$}
\end{picture}
\caption{Grid points and discrete values for bi-domain equations}
\label{bid1}
\end{figure}

We introduce a fixed time step $\Delta t>0$ and define discrete times $t_n=n\Delta t$ for $n\in\N_0$.
Thus we obtain nodal grid point $(x_j,t_n)\in [0,1]\times\R_{\ge 0}$ for our computational domain.
For the functions $u$ and $v$ introduced in the previous section, we set nodal values to be
$$
u_{j}^n \approx u(x_{j},t_n)\quad \text{for}\; j=0,1,...,m-1,
\quad v_{j}^n\approx v(x_{j},t_n) \quad\text{for}\;j=m+1,m+2,...,N
$$
as well as at the interface node $u_{m}^n\approx u(x_{m_{-}},t_n)$ and 
$v_{m}^n\approx v(x_{m_{+}},t_n)$. 

The discretization of the diffusion equation is achieved via an explicit forward time step of length $\Delta t>0$ and 
the standard central difference in space for the second derivative in space. This gives us the forward in time central in space
(FTCS) scheme. 

We introduce the parameters $\nu_-=D_{-}\frac{\Delta t}{(\Delta x)^2}$ and $\nu_+=D_{+}\frac{\Delta t}{(\Delta x)^2}$.
The updates determining the interior node values for the system \eqref{eqn1} written as explicit iterations
are then given as
\begin{eqnarray}
\label{scheme_u}
u^{n+1}_{j}&=&u_j^n +\nu_- (u^n_{j+1}-u_{j}^{n})-\nu_-(u_j^n- u_{j-1}^{n}) \qquad\mbox{for}\;j=1,...,m-1\\
\label{scheme_v}
v^{n+1}_{j}&=&v_j^n +\nu_+ (v^n_{j+1}-v_{j}^{n}) -\nu_+ (v_j^n-v_{j-1}^{n}) \qquad\mbox{for}\;j=m+1,...,N-1.
\end{eqnarray}
It is well-known by 
stability analysis, see e.g.\ Morton and Mayers \cite{bMOMA} or Thomas \cite{bTHO}, 
that the time steps of this scheme are restricted by 
$\Delta t\le \frac{(\Delta x)^2}{2 \max\{D_{-},D_{+}\}}$ giving $0<\nu_\pm< \frac 12$.
Due to the square term, the time steps become very small on fine meshes. This leads to the preference for implicit schemes in practice.

\subsection{Numerical homogeneous Neumann boundary conditions}

In order to determine $u_0^{n+1}$ and $v_N^{n+1}$ we need numerical boundary conditions. We cannot use the
formulas of type \eqref{scheme_u} for $j=0$ or \eqref{scheme_v} for $j=N$ directly. 
To compute these nodal values we introduce the ghost points $x_{-1} =-\Delta x$ and $x_{N+1} =(N+1)\Delta x$ as well as 
the corresponding ghost values $u_{-1}^n$ and $v_{N+1}^n$.
We implement discrete homogeneous outer Neumann boundary conditions. 

Let us consider the boundary at $x=0$. The simplest
two finite difference approximations we can use, namely the one-sided and the central difference with respect
to the node $x_0$, are
\begin{equation}
\label{bdry_fd}
\frac{u_0^n-u_{-1}^n}{\Delta x} = 0 \qquad \mbox{and}\qquad \frac{u_1^n-u_{-1}^n}{2\Delta x}=0.
\end{equation}
These give the values $u_{-1}^n=u_0^n$ and $u_{-1}^n = u_1^n$ respectively. They are inserted into the
updates \eqref{scheme_u} for $j=0$. Analogously, we proceed for $j=N$ using \eqref{scheme_v}. The one-sided differences give
\begin{equation}
\label{onesided_bdry}
u^{n+1}_{0}=u_{0}^{n}+\nu_- (u^n_{1}-u_{0}^{n}) \qquad\mbox{and}\qquad v^{n+1}_{N}=v_{N}^{n}-\nu_+ (v^n_{N}-v_{N-1}^{n}).
\end{equation}
The central differences lead to extra factors of $2$ in the differences.
\begin{equation}
\label{central_bdry}
u^{n+1}_{0}=u_{0}^{n}+2\nu_- (u^n_{1}-u_{0}^{n})\qquad\mbox{or}\qquad v^{n+1}_{N}=v_{N}^{n}-2\nu_+( v^n_{N}-v_{N-1}^{n}).
\end{equation}
We will see in Section \ref{sec:cons} that the conservation property will determine when 
each of the formulas \eqref{onesided_bdry} or \eqref{central_bdry}
is useful. We will see for the nodal-based scheme that only \eqref{central_bdry} correctly represents
the conservation property. 
Later, we will also consider a finite volume type scheme. It needs \eqref{onesided_bdry}
as boundary conditions. In finite volume schemes, the one-sided difference formula turns out to be a central 
difference with respect to the cell
boundary at $x=0$ or $x=1$. This gives the correct boundary fluxes there. 

In the literature, these boundary conditions seem to have only been compared
in terms of the fact that \eqref{central_bdry} is a second-order approximation in space and \eqref{onesided_bdry}
only of first order, while the interior updates \eqref{scheme_u} and \eqref{scheme_v} are second order in space, 
see e.g.\ Hundsdorfer and Verwer \cite[Subsection I.5.3]{bHUVE} 
or Thomas \cite[Section 1.4]{bTHO}. In \cite{bHUVE}, there is
a nice symmetry argument for using \eqref{central_bdry}. For $L^p$ error estimates, the truncation error of the boundary
discretization is allowed to be one order lower than in the interior. The total estimates are determined by summing up
the cells or elements. The number of those that contain boundary conditions is one order of $\frac 1h$ less than those
in the interior. Therefore, the overall order in the space of a scheme using \eqref{onesided_bdry} 
is maintained despite the lower
order of the truncation error. So, the order of the scheme does not make much of a distinction 
between the two numerical boundary conditions.

For a complete scheme, it remains to determine updates $u_{m}^{n+1}$ and $v_m^{n+1}$ 
via various numerical coupling conditions below.
We will obtain the numerical coupling conditions results analogous to the numerical boundary conditions for the 
discrete numerical coupling conditions via two fluxes.
The Dirichlet-Neumann coupling differs due to the Dirichlet condition. 

\subsubsection*{The piecewise linear finite element method}

The above updates \eqref{scheme_u}, \eqref{scheme_v} and \eqref{central_bdry} can be achieved on our regular mesh by 
taking piecewise linear finite elements as the spatial semi-discretization. The finite element functions are
linear on the intervals $[x_{j-1},x_j]$ for $j=1,\ldots ,N$, continuous at the nodes and with nodal values as above.
In the case of flux couplings below, we would allow a discontinuity at $x_m$. But let us ignore that for the moment.
After quadrature, this finite element method gives a system of ordinary
differential equations in time. Let $\ul w(t)\in\R^{N+1}$ be the vector of nodal values of the piecewise linear
approximations at time $t\ge 0$. Then the resulting system has the form $\ul M\dot{\ul w}(t) +\ul A\ul w(t)=\ul 0$ with
the mass matrix $\ul M\in\R^{(N+1)\times (N+1)}$ and the stiffness matrix $\ul A\in\R^{(N+1)\times (N+1)}$ 
that are calculated using the hat basis functions.
In the literature, these tri-diagonal matrices can be found as $(N-1)\times(N-1)$ matrices 
for the single domain heat equation 
with zero Dirichlet boundary data. The mass matrix has the entries $\Delta x\frac 46$ on the diagonal and
$\Delta x\frac 16$ on the two secondary diagonals, analogously the entries are $\frac 2{\Delta x}$ 
and $\frac {-1}{\Delta x}$ for the stiffness matrix, see e.g.\ Wait and Mitchell \cite[Section 5.2]{bWAMI}.

The homogenous Neumann condition is automatically satisfied for finite elements. 
This comes from the variational principle behind them.
Dirichlet conditions need to be enforced by putting boundary values to zero. A flux boundary condition is
put into the underlying functional for the Galerkin method. In case the boundary fluxes vanish, nothing has to be done. 
For this reason, Neumann boundary conditions
are also called natural or do nothing boundary conditions, see e.g.\ Johnson \cite[Section 1.7]{bJOHN}.
The homogeneous Neumann condition requires the use of the nodal values at $x_0$ and $x_N$. The mass and stiffness matrices
are then tri-diagonal $(N+1)\times (N+1)$ matrices. The entries for the interior nodes $x_j$ for $j=1,\ldots, N-1$ are the same
as those for $j=2,\ldots ,N-2$ in the case of the Dirichlet boundary condition. But for $j=0,N$ the homogeneous Neumann boundary
condition gives the entries $\Delta x\frac 26$ on the diagonal and
$\Delta x\frac 16$ on the two secondary diagonals, analogously the entries are $\frac 1{\Delta x}$ 
and $\frac {-1}{\Delta x}$ for the stiffness matrix. This comes from the fact that only half of a hat basis function is
used at the boundary nodes for the Neumann boundary data.

In computations, the mass matrix $M$ is dealt with using a commonly applied procedure called mass lumping, 
see Hundsdorfer and Verwer \cite[Section III.5]{bHUVE}, Quarteroni and Valli \cite[Section 11.4]{bQUVA1}
or Thom\'ee \cite[Chapter 15]{bTHOME}. 
One adds in each row all entries of $\ul M$ and puts them as entries into the regular diagonal
matrix $\ul D$. Thus, the need to solve a system of linear equations in each time step is eliminated. 
This gives the system $\dot{\ul w}(t) +\ul D^{-1}\ul A\ul w(t)=\ul 0$. An exact quadrature to obtain $\ul M$
requires the Simpson rule on each interval since the integrands are quadratic. Taking the trapezoidal rule as an
approximate quadrature gives the lumped matrix. Quarteroni and Valli, as well as Thom\'ee
discussed this for the case of two space dimensions. Thom\'ee gave a second
interpretation of the procedure.
The matrix $\ul D$ for the Neumann boundary conditions has the entries $1$ for the interior nodes and $\frac 12$ for the boundary
nodes. Then a forward time step leads to the boundary updates \eqref{central_bdry} with the factor of $2$ 
as well as the interior updates \eqref{scheme_u}, \eqref{scheme_v}.
We will see in Section \ref{sec:cons} that the finite element method 
is consistent with the conservation property.

\subsection{Discrete Dirichlet-Neumann coupling conditions}

As a first step we consider the Dirichlet-Neumann conditions \eqref{nnm} for the bi-domain diffusion model, i.e.\ the Dirichlet condition
$u(\frac 12,t)=v(\frac 12,t)$ and the Neumann condition 
$D_{-}\frac{\partial  u(\frac 12,t)}{\partial x} =D_{+}\frac{\partial  v(\frac 12,t)}{\partial x}$.
We must choose one of them in order to determine $u_m^{n+1}$ and the other for $v_m^{n+1}$.
We take the Neumann condition for $u_m^{n+1}$ as our first step and discretize it using forward differences and a ghost point value $u_{m+1}^n$ as
$D_{-}\frac{u_{m+1}^n-u_{m}^n}{\Delta x}=D_{+}\frac{v^n_{m+1}-v^n_{m}}{\Delta x}$. 
This is solved for $u_{m+1}^n-u_{m}^n$ and the result inserted into \eqref{scheme_u} for $j=m$. Using $\nu_+=\frac{D_{+}}{D_{-}}\nu_-$
and $u_m^n=v_m^n$ we obtain the following update
\begin{equation}
u^{n+1}_{m}= u^n_{m}+ \nu_{+} (v_{m+1}^n-u_{m}^n)-\nu_{-}(u_{m}^n-u_{m-1}^n).
 \label{f16}
\end{equation}
Note that Giles \cite{GIL} suggested a derivation using the fluxes to arrive at the same update.

Then we use the Dirichlet condition to set $v_m^{n+1}=u_m^{n+1}$ as our second step to achieve the coupling.
We assume that the discrete initial data satisfy $v_m^0=u_m^0$. So, we will always have $v_m^n=u_m^n$, and
we could eliminate the variable $v_m^n$ from explicit schemes using the Dirichlet coupling condition.
We keep it for comparison with conditions where this is not the case.

Using the backward differences above means that one introduces the ghost value $v_{m-1}^n$. This can be inserted
into \eqref{scheme_v} and with the Dirichlet condition $u_m^n=v_m^n$ gives the same scheme, the Dirichlet update
becoming $u_m^{n+1}=v_m^{n+1}$. Taking a central difference here does not make much sense, since we would be introducing
two ghost values into one Neumann condition.

Now, the fully discretized explicit FTCS nodal scheme for the bi-domain diffusion model with Dirichlet-Neumann
coupling conditions is given as 
\begin{eqnarray}
\label{scheme_c1}
u_{0}^{n+1}=&u_{0}^n+2\nu_- (u_{1}^{n}-u_{0}^n) \quad\qquad\qquad\qquad\qquad&\quad\mbox{for $j=0$},\nonumber \\
u_{j}^{n+1}=&u_j^n +\nu_-( u_{j+1}^{n}-u^n_j) -\nu_-(u_j^n-u^n_{j-1}) \;\quad&\quad\mbox{for $0<j<m$},\nonumber \\
v_{j}^{n+1}=&v_j^n+\nu_+ (v_{j+1}^{n}-v^{n}_j)-\nu_+(v_j^n-v^n_{j-1})\;\quad &\quad\mbox{for $m<j<N$}, \nonumber \\
v_{N}^{n+1}=&v_{N}^n-2\nu_+ (v_{N}^{n}-v_{N-1}^n) \quad\qquad\qquad\qquad&\quad\mbox{for $j=N$},
\end{eqnarray}
with the coupling conditions
\begin{eqnarray}
\label{dn_disc}
u_{m}^{n+1}=& u^n_{m}+\nu_{+}(v^n_{m+1}-u^n_{m})-\nu_- (u^n_{m}-u^n_{m-1})& \quad \mbox{for $j=m$},\nonumber \\
v^{n+1}_m=&u^{n+1}_m  \quad\qquad\qquad\qquad\qquad\qquad\qquad\qquad&\quad\mbox{for $j=m$}.\nonumber \\
\end{eqnarray}

\subsubsection*{The coupling scheme of Giles}

Giles \cite{GIL} considered heat diffusion with a heat capacity $c$ and conductivity $k$, i.e.\ in our notation an
equation of the form $cu_t-ku_{xx}=0$. He further considered 
different mesh sizes $\Delta x_{\pm}$ in the two sub-domains. 
Note that Giles \cite[Eq.~(26)]{GIL} did not distinguish $u$ and $v$ in the bi-domain case, and he did not
introduce double values at the interface node. Therefore, 
$v^{n+1}_m=u^{n+1}_m$ is implicitly automatically implied in his 
scheme. As we mentioned, this is something we could have also done above without changing the outcome of the explicit computations. 
But, we will need our approach with double values at the interface for the other 
coupling conditions that do not include the Dirichlet condition.

Before making a supposed simplification, Giles had derived a correct coupling scheme \cite[Eqs.~(15),(16)]{GIL}. 
In our notation, the update is given as
\begin{equation}
\label{giless2}
u_{m}^{n+1}= u^n_{m}+\frac{2r\nu_{+}}{1+r}(v^n_{m+1}-u^n_{m})-\frac{2\nu_{-}}{1+r}(u^n_{m}-u^n_{m-1})
\end{equation}
where $r=\frac{c_-\Delta x_-}{c_+\Delta x_+}$ and $\nu_{\pm}=\frac{k_{\pm}\Delta t}{c_{\pm}(\Delta x_{\pm})^2}$. 
Note that setting $c_{\pm}=1$, $k_{\pm}=D_{\pm}$, $\Delta x_-=\Delta x_+$ and $r=1$ in \eqref{giless2} gives our formula \eqref{f16}.
 
However, in Gile's stability analysis, \cite[(26)]{GIL} used the update for
\begin{equation}
\label{giless1}
u_{m}^{n+1}= u^n_{m}+2 r\nu_{+}(v^n_{m+1}-u^n_{m})-2\nu_- (u^n_{m}-u^n_{m-1}).
\end{equation}
This discretization was obtained by introducing an inconsistency in the time discretization. It leads to the loss of 
the conservation property, see Section \ref{sec:cons}, and some instabilities, see Giles \cite{GIL}.

\subsection{Discrete other coupling conditions}

We discretize the heat flux coupling conditions defined in \eqref{eq111} via an explicit discretization method with one-sided 
differences. The heat flux coupling conditions are 
$D_{-}\frac{\partial u}{\partial x}= D_{+}\frac{\partial v}{\partial x}=H(v-u).$
For the nodes $j\ne m$, we use the formulas in \eqref{scheme_c1}. Only the coupling conditions for $j=m$ will be replaced.
Now for the interface node $j=m$, the updates \eqref{scheme_u} and \eqref{scheme_v} are
\begin{eqnarray}
 \label{n32}
 u^{n+1}_{m}&=&u^{n}_{m}+\nu_-(u^{n}_{m+1}-u^{n}_{m})-\nu_- (u^{n}_{m}-u^n_{m-1})\nonumber\\
 v^{n+1}_{m}&=&v^{n}_{m}+\nu_+(v^{n}_{m+1}-v^{n}_{m})-\nu_+(v^{n}_{m}-v^n_{m-1}).
\end{eqnarray}
In these two formulas, we have two ghost point values $u^n_{m+1}$ and $v^n_{m-1}$. We calculate these values by discretizing
the two coupling conditions \eqref{eq111}.
For the first one, we take the forward difference approximation, and for
the second one is the backward difference approximation to obtain these ghost point values
\begin{equation}
D_{-}\frac{u_{m+1}^n-u_{m}^n}{\Delta x}=D_{+}\frac{v_{m}^n-v_{m-1}^n}{\Delta x}= -J_{heat}(u_m^n,v_m^n) = H(v_{m}^n-u_{m}^n).
\label{y7}
\end{equation}
Solving the first equation gives $u_{m+1}^{n}=u_{m}^n+\frac{H\Delta x}{D_{-}}(v_{m}^n-u_{m}^n)$.
We substitute into \eqref{n32} and obtain using $\nu_-\frac{H\Delta x}{D_{-}} 
=\frac{ D_{-}\Delta t}{H(\Delta x)^2}\frac{\Delta x}{D_{-}} =\frac{H \Delta t}{\Delta x}$ the update
\begin{eqnarray} 
u_{m}^{n+1}&=&u^n_m+\nu_- \Big(u_{m}^n+\frac{H\Delta x}{D_{-}}(v_{m}^n-u_{m}^n)-u^n_m\Big)-\nu_-( u_{m}^{n}-u^n_{m-1})\nonumber\\
&=&u^n_m-\nu_-(u^n_m-u^n_{m-1})+\frac{H \Delta t}{\Delta x}(v^n_m-u^n_m),\nonumber\\
&=&u^n_m-\nu_-(u^n_m-u^n_{m-1})-\frac{\Delta t}{\Delta x}J_{heat}(u^n_m,v_m^n).
\label{300}
\end{eqnarray}
Solving the second equation of \eqref{y7} we get $v_{m-1}^{n}=v_{m}^n-\frac{H \Delta x}{D_{+}} (v_{m}^n-u_{m}^n)$ and
analogously
\begin{equation}
v_{m}^{n+1}=v^n_m+\nu_+(v^n_{m+1}-v^n_m)+\frac{\Delta t}{\Delta x}J_{heat}(u^n_m,v_m^n).
\label{3001}
\end{equation}
For the coupled scheme, we proceed analogously as in \eqref{scheme_c1}. Only we are replacing the 
updates for the Dirichlet-Neumann coupling at $j=m$
by the new formulas  \eqref{300} and \eqref{3001} for the heat flux coupling conditions.

Now, we want to discretize coupling conditions \eqref{eq111} via central difference approximations as
\begin{equation}
 D_{-}\frac{u_{m+1}^{n}-u_{m-1}^n}{2\Delta x}=D_{+}\frac{v_{m+1}^{n}-v_{m-1}^n}{2\Delta x}=H(v^n_m-u^n_m).
\end{equation}
This gives $u_{m+1}^{n}=u_{m-1}^n+\frac{2H\Delta x}{D_{-}}(v_{m}^n-u_{m}^n)$.
Substituting into \eqref{n32} gives additional factors of $2$ in the updates
\begin{eqnarray}
\label{2001}
u_{m}^{n+1}&=&u^n_m-2\nu_-(u^n_{m}-u_{m-1}^n)-\frac{2\Delta t}{\Delta x} J_{heat}(u^n_m,v_m^n)\nonumber\\
v_{m}^{n+1}&=&v^n_m+2\nu_+ (v_{m+1}^n-v^n_m)+\frac{2\Delta t}{\Delta x}J_{heat}(u^n_m,v_m^n).
\end{eqnarray}
Note the analogy to the numerical homogeneous Neumann conditions \eqref{central_bdry}. These are the coupling
conditions that will prove to maintain the conservation property for the scheme \eqref{scheme_c1}. The coupling conditions
\eqref{2001} replace \eqref{dn_disc}.

We now refer to the other coupling conditions defined in Subsection \ref{subsec:couple}.  
We obtain, for example, the explicit discretization for the channel pumping conditions
by replacing the heat flux $J_{heat}(u_m^n,v_m^n)=-H(v^n_m-u^n_m)$ 
by the channel flux $J_{ch}(u_m^n,v_m^n)=\Psi \frac{u^n_m-\alpha v^n_m}{\beta+\gamma u^n_m+\delta v^n_m}$ 
in the updates \eqref{300} and \eqref{3001} or in \eqref{2001}. Analogously, we proceed to the membrane coupling conditions.

\subsection{An explicit finite volume type scheme}
\label{sub_fv}

\begin{figure}[htp]
\setlength{\unitlength}{1cm}
\begin{picture}(15,2) 
\put(2.0,1){\line(1,0){12}} 
\put(1.9,0.3){$0$}
\put(2.0,0.7){\line(0,1){0.6}}    
\put(2.7,1){\circle*{0.1}}
\put(2.5,0.5){$x_1$}
\put(2.5,-0.1){$\sigma_1$}
\put(2.4,1.5){$u_1^n$}
\put(3.4,0.7){\line(0,1){0.6}}
\put(4.1,1){\circle*{0.1}}
\put(3.9,0.5){$x_2$}
\put(3.9,-0.1){$\sigma_2$}
\put(4,1.5){$u_2^n$}
\put(4.8,0.7){\line(0,1){0.6}}
\put(5.7,0.5){$\cdots$}
\put(6.6,0.7){\line(0,1){0.6}}
\put(7.3,1){\circle*{0.1}}
\put(7.1,0.5){$x_m$}
\put(7.1,-0.1){$\sigma_m$}
\put(7.1,1.5){$u_m^n$}
\put(8,0.7){\line(0,1){0.6}}
\put(7.9,0.2){$\frac 12$}
\put(8.7,1){\circle*{0.1}}
\put(8.5,0.5){$x_{m+1}$}
\put(8.5,-0.1){$\sigma_{m+1}$}
\put(8.5,1.5){$v_{m+1}^n$}
\put(9.4,0.7){\line(0,1){0.6}}
\put(10.8,0.5){$\cdots$}
\put(12.6,0.7){\line(0,1){0.6}}
\put(13.3,1){\circle*{0.1}}
\put(13.1,0.5){$x_N$}
\put(13.1,-0.1){$\sigma_N$}
\put(13.1,1.5){$v_N^n$}
\put(13.9,0.3){$1$}
\put(14.0,0.7){\line(0,1){0.6}} 
\end{picture}
\caption{Cells, nodes, and cell or nodal values for the finite volume scheme.}
\label{fbb1}
\end{figure}
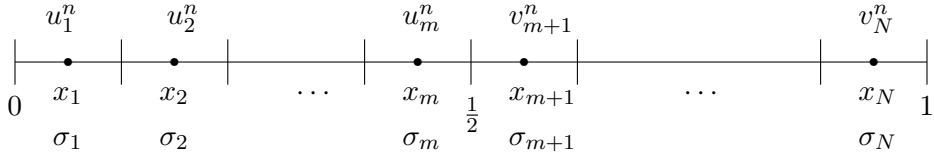

We now want to consider a finite volume type scheme. These types of schemes are quite popular for compressible fluid flow computations
and useful when the conservation property of quantities comes into play. 
For this scheme we define for $m\in\N$ the $N=2m$ cells of length $\Delta x=1/N$ with midpoints $x_j = \Delta x(j-1/2)$ and 
boundary points $x_{j\pm\frac{1}{2}}=x_j\pm \frac{\Delta x}{2}$
for $j=1,\ldots ,N$. The cells are the sub-intervals $\sigma_j=[x_{j-\frac{1}{2}}, x_{j+\frac{1}{2}}]$, see
Figure \ref{fbb1}. The number of cells and nodes is always even in our coupling problems. Each sub-domain $[0,\frac 12]$
and $[\frac 12,1]$ has $m$ cells. The nodal points of
the nodal-based schemes considered above
have become the cell boundary points.  

On the cells $\sigma_j$, we consider the solutions to be constant and assign the values $u_j^n$ for $j=1,\ldots ,m$
as wells as $v_j^n$ for $j=m+1,\ldots ,N$ to the nodes $x_j$ at the cell center. Since the interface sits at cell boundaries,
we do not have to use a node with a double value for the coupling.

We seek approximations of $u$ and $v$ by integral averages on the cells to represent our solutions, i.e.\ by
$u_{j}^n\approx \frac{1}{\Delta x}\int_{x_{j-\frac{1}{2}}}^{x_{j+\frac{1}{2}}} u(x,t_n)\,dx$ for $j=1,...,m$
and $v_{j}^n\approx\frac{1}{\Delta x}\int_{x_{j-\frac{1}{2}}}^{x_{j+\frac{1}{2}}} v(x,t_n)\,dx$ 
for $j=m+1,...,N=2m$.
For the discretization of the initial data, we can use these integral averages. The integrals may be
replaced by a quadrature rule.

Finite volume schemes use a numerical flux formulation that is very useful in the presence of the conservation property.
Using it, they maintain this property automatically. We introduce the right and left-hand numerical flux functions for 
cell $\sigma_j$ for $j=1,\ldots,m$ as
\begin{equation}
\label{num_flux}
F_\pm^{j+\frac 12,n}(u_{j+1}^n,u_j^n) = \nu_\pm(u_{j+1}^n-u_j^n)\qquad\text{and} 
\qquad F_\pm^{j-\frac 12,n}(u_j^n,u_{j-1}^n) = \nu_\pm(u_j^n-u_{j-1}^n).
\end{equation}
The numerical fluxes for $j=m+1,\ldots, N$ are analogous with $u_j^n$ replaced by $v_j^n$.
Note that usually $\frac{\Delta t}{(\Delta x)^2}$ is not included in the fluxes for the purpose of numerical analysis. 
But we find it convenient to do so here. An important property is that $F_\pm^{j+\frac 12,n}(u_{j+1}^n,u_j^n) =
F_\pm^{(j+1)-\frac 12,n}(u_{(j+1)}^n,u_{(j+1)-1}^n)$, i.e.\ the right hand flux for cell $\sigma_j$ is equal to the left hand
flux of cell $\sigma_{j+1}$. In the updates, they will appear with opposite signs. 
With these numerical flux functions, we obtain the flux form of the updates. 
\begin{eqnarray*}
u_j^{n+1} &=& u_j^n + F_-^{j+\frac 12,n}(u_{j+1}^n,u_j^n) -F_-^{j-\frac 12,n}(u_j^n,u_{j-1}^n) \quad\text{for}\; j=2,\ldots ,m-1,\\
v_j^{n+1} &=& u_j^n + F_+^{j+\frac 12,n}(v_{j+1}^n,v_j^n) -F_+^{j-\frac 12,n}(v_j^n,v_{j-1}^n) \quad\text{for}\; j=m+2,\ldots ,N-1=2m-1.
\end{eqnarray*}
For the interior cells $\sigma_j$ for $1<j<m$ and $m+1<j<N$, this gives the FTCS updates used in \eqref{scheme_u} and \eqref{scheme_v}. 

As numerical boundary conditions for the boundary cells $\sigma_1$ and $\sigma_N$,
we use the one-sided difference formulas from \eqref{onesided_bdry}. For the left boundary, we have to replace $j=0$ by $j=1$, e.g.\ 
$u_1^{n+1} =u_1^n +F_-^{\frac 32,n}(u_2^n,u_1^n)$. This will be justified in Section \ref{sec:cons}. 

In finite difference form, the explicit FTCS finite volume type scheme is
\begin{eqnarray}
\label{fvol2}
u^{n+1}_{1}&=&u^n_1+\nu_-( u^n_{2}-u_{1}^{n})\nonumber\\
u^{n+1}_{j}&=&u_j^n+\nu_-( u^n_{j+1}-u_{j}^{n}) -\nu_-(u_j^n- u_{j-1}^{n}) \qquad\mbox{for}\;j=2,...,m-1\nonumber\\
v^{n+1}_{j}&=&v_j^n+\nu_+( v^n_{j+1}-v_{j}^{n}) - \nu_+ (v_j^n-v_{j-1}^{n} )\qquad\mbox{for}\;j=m+2,...,N-1\\
v^{n+1}_{N}&=&v_N^n-\nu_+( v^n_{N}-v_{N-1}^{n}).\nonumber
\end{eqnarray}

Now we derive the discretization scheme of the Dirichlet-Neumann coupling \eqref{nnm} for our finite volume type scheme using
ghost cell values $v^n_m$ and $u^n_{m+1}$.
We discretize the Neumann coupling condition via the central 
differences with respect to the boundary. These are one-sided differences for $u_m^n$ and $v_{m+1}^n$. We set
$D_{+}\frac{v^n_{m+1}-v^n_{m}}{\Delta x}=D_{-}\frac{u^n_{m+1}-u^n_m}{\Delta x}$.
Using the Dirichlet condition $u^n_m=v^n_m$  
this implies that $\nu_{-}(u^n_{m+1}-u^n_m)=\nu_{+}(v^n_{m+1}-v^n_m)=\nu_{+}(v^n_{m+1}-u^n_m)$. 
This gives for $j=m, m+1$ the Dirichlet-Neumann updates
\begin{eqnarray*}
u^{n+1}_{m}&=&u^n_m+\nu_{-}(u^n_{m+1}-u^n_m)-\nu_{-}(u^n_m-u^n_{m-1})\\ 
&=&u^n_m+\nu_+(v^n_{m+1}-u^n_m)-\nu_-(u^n_m-u^n_{m-1}),\\
v^{n+1}_{m+1}&=&v_{m+1}^n+\nu_+ (v^n_{m+2}-v_{m+1}^{n}) - \nu_+ (v_{m+1}^n-u_m^{n} )
\end{eqnarray*}
Note that the first formula for $u_m^n$ is the same as in the nodal-based case \eqref{scheme_c1}. 
A difference is only in the outer boundary conditions and the interpretation of the discrete values.

For the various flux coupling conditions, we take the updates \eqref{300} and \eqref{3001} and insert the appropriate fluxes.
The justification will be given in the next section. The numerical flux functions for $j=m+\frac 12$ 
then contain the added extra term $\frac {\Delta t}{\Delta x}J(u_m^n,v_m^n)$.
\section{The discrete conservation property for bi-domain diffusion equations}
\label{sec:cons}

Consider our coupling problem \eqref{eqn1} with homogeneous Neumann boundary data.
The total amount of concentration at any time $t$ is $C(t) =\int_0^{\frac 12}u(s,t)\d s+\int_{\frac 12}^1v(s,t)\d s$.
Differentiating, then using the equations \eqref{eqn1} and the boundary conditions, we obtain
\begin{eqnarray*}
\frac \d{\d t} C(t) &=&\int_0^{\frac 12}\frac\partial{\partial t}u(s,t)\d s+\int_{\frac 12}^1\frac\partial{\partial t}v(s,t)\d s
=D_-\int_0^{\frac 12}\frac{\partial^2}{\partial s^2}u(s,t)\d s+D_+\int_{\frac 12}^1\frac{\partial^2}{\partial s^2}v(s,t)\d s\\
&=& D_- \frac{\partial}{\partial x}u(\frac 12,t) -D_+\frac{\partial}{\partial x}v(\frac 12,t) =0
\end{eqnarray*}
if the coupling condition $D_- \frac{\partial}{\partial x}u(\frac 12,t) =D_+\frac{\partial}{\partial x}v(\frac 12,t)$ holds.
This is true for all the conditions that we use. Note that the Dirichlet condition or the nature of fluxes is not relevant.
The result means that $C(t)$ remains constant in time. Obviously, this property should be shared by the discrete schemes that we use.
It should hold up to machine accuracy in any computation. Behind this conservation property are important physical principles such as
mass conservation for concentrations or energy conservation in the case of the heat equation.

\subsection{Discrete mass conservation}

In order to check whether the conservation property holds throughout the computational time interval, we will
compute the discrete total mass $C_{total}$ in the domain at each time step $t_n=n\Delta t$ for $n\in\N_0$.

\subsubsection*{The nodal-based scheme}

For this purpose, we have to introduce cells for our nodal-based scheme. The nodal value then represents a constant
value over the whole cell. The total amount of our discrete variables on a cell is the value times the length of the cell. 
We define cell boundary points $x_{j\pm\frac 12} =j\Delta x\pm\frac{\Delta x}2$ for any $j=1,\cdots , N-1$
and introduce the cells $\sigma_j =[x_{j-\frac 12},x_{j+\frac 12}]$
for $j=1,\cdots ,m-1,m+1,\cdots N-1$. This is analogous to our finite volume type scheme. But the nodes and cells 
are shifted by $\frac{\Delta x}2$. For the boundary nodes
we take the cells to be $\sigma_0 =[x_0,x_{0+\frac 12}]$ and $\sigma_N=[x_{N-\frac 12}, x_N]$. At the interface, we introduce
two cells $\sigma_{m-}=[x_{m-\frac 12},x_{m-}]$ and $\sigma_{m+}=[x_{m+},x_{m+\frac 12}]$. 
These smaller cells are the reason why
the flux boundary and coupling conditions are different for the schemes.

For $n\in \mathbb{N}_0$, we define concentration sum
\begin{equation}
\label{dnc2}
C_n=\frac{1}{2}u^n_0+u^n_1+\cdots +u^n_{m-1}+\frac{1}{2}u^n_m+\frac{1}{2}v^n_m+v^n_{m+1}
+\cdots+v^n_{N-1}+\frac{1}{2}v^n_{N}
\end{equation}
and the discrete total concentrations $\overline{C}_n = \Delta x C_n$.
The integral of the initial concentration over the whole domain gives the initial total concentration $\overline{C}_0$.
If we take the initial nodal values to be the integral averages over the cells of the exact initial data, 
then they give the same initial total concentration $\overline{C}_0 = \Delta x C_0$. 
Any approximation of the integral averages, e.g.,  by quadrature, will introduce a small initial error in total mass.
The common factor $\Delta x$ is not needed for the considerations that follow. Therefore, we will work with $C_n$.
The { discrete conservation property} for a nodal-based scheme then is that $C_{n+1} = C_n$ for $n\in\N_0$.
An attempt to introduce the discrete conservation property was made in Morton and Mayers \cite[Section 2.14]{bMOMA}, 
but not quite in an adequate way and not put to much use.

\noindent
{\bf Boundary conditions:}
It is easily seen that the numerical homogeneous Neumann boundary conditions \eqref{central_bdry} obtained via central 
differencing give the discrete conservation property. The same is true for Neumann boundary conditions 
with non-zero fluxes. The boundary conditions \eqref{onesided_bdry} using the 
one-sided differences produce
errors of the order $\frac{\nu_-} 2(u_1^n-u_0^n)$ and $\frac{\nu_+} 2(v_N^n-v_{N-1}^n)$.

\noindent
{\bf Dirichlet-Neumann coupling:}
We will need the above form of discretization in which the cell $\sigma_m = [x_{m-\frac{1}{2}}, x_{m+\frac{1}{2}}]$ 
is split into two sub-cells for the
coupling conditions that do not involve the Dirichlet condition.
When we consider the discrete Dirichlet condition $u_m^n =v_m^n$, 
the term $\frac{1}{2}u^n_m+\frac{1}{2}v^n_m$ is replaced by $u_m^n$.
This is the reason why the Dirichlet-Neumann coupling does not need the extra factors 
of $2$ that are present in the boundary conditions.

We will make use of the numerical flux functions \eqref{num_flux}.
With these, we can rewrite the nodal-based schemes. 
For instance \eqref{scheme_u} and the boundary condition \eqref{central_bdry} for $j=0$ are
$$
u_j^{n+1} = u_j^n + F_-^{j+\frac 12,n}(u_{j+1}^n,u_j^n)-F_-^{j-\frac 12,n}(u_j^n,u_{j-1}^n)
\qquad\text{and}\qquad u_0^{n+1} =u_0^n + 2F_-^{\frac 12,n}(u_{1}^n,u_0^n)
$$
Analogously, we can treat the other updates and use $F_-^{j+\frac 12,n}(v_{j+1}^n,v_j^n)$ as well as
$F_-^{j-\frac 12,n}(v_j^n,v_{j-1}^n)$ for $j\ge m$. This is basically a finite volume type formulation
of the nodal-based scheme using shifted cells.
Since $F_\pm^{j+\frac 12,n}(u_{j+1}^n,u_j^n) = F_\pm^{(j+1)-\frac 12,n}(u_{(j+1)}^n,u_{(j+1)-1}^n)$ 
and $F_\pm^{j-\frac 12,n}(u_j^n,u_{j-1}^n)=F_\pm^{(j-1)+\frac 12,}(u_{(j-1)+1}^n,u_{(j-1)}^n)$ we see that the
fluxes cancel when the updates are inserted into $C_{n+1}$ and summed up. The updates for $j=0,m,N$ need the extra factor of $2$
in the flux to compensate for the half-sized cells there. Then we obtain $C_{n+1} =C_n$. This is the discrete conservation
property.

For the Dirichlet-Neumann coupling, we use the cell $\sigma_m$ and $\frac{1}{2}u^n_m+\frac{1}{2}v^n_m =u_m^n$ in \eqref{dnc2}.
We see that \eqref{scheme_c1} can be written just using the fluxes
\eqref{num_flux} and all fluxes cancel to give $C_{n+1} =C_n$.
Giles \cite{GIL} used \eqref{giless1} which for $r=1$ contains an extra factor of $2$ and therefore gives
$C_{n+1} = C_n -\nu_-(u_m^n-u_{m-1}^n) +\nu_+(v_{m+1}^n-u_m^n)$. So, it does not have the conservation property.
The correct scheme \eqref{giless2} is identical to \eqref{scheme_c1} for $r=1$. So, it has the conservation property
in this case. For the case $r\ne 1$ in \eqref{giless2}, one would have to modify \eqref{dnc2} to take the different mesh sizes 
into account. This will be done in a further paper \cite{CMW3}.

\noindent
{\bf Flux couplings:}
For the heat flux coupling, we have to take the central difference coupling \eqref{2001} to compensate for $\frac 12u_m^n$ and
$\frac 12v_m^n$ in $C_{n+1}$. We have
\begin{eqnarray*}
u_m^{n+1} &=& u_m^n -2 F_-^{m-\frac12}(u_m^n,u_{m-1}^n) -\frac{2\Delta t}{\Delta x} J_{heat}(u^n_m,v_m^n),\\
v_m^{n+1} &=& v_m^n +2 F_-^{m+\frac12}(v_{m+1}^n,v_m^n) +\frac{2\Delta t}{\Delta x} J_{heat}(u^n_m,v_m^n).
\end{eqnarray*}
The extra heat flux terms cancel in the summation in $C_{n+1}$ giving $C_{n+1}=C_n$. 
Obviously, this will hold for all other discrete coupling conditions.

\subsubsection*{Piecewise linear finite elements}

Note that for piecewise linear finite elements on our nodal mesh, one can compute the total integral
of the solution over the interval $[0,1]$ using exact quadrature by the mid-point or trapezoidal rules as
\begin{eqnarray*}
\overline{C}_n&=&\Delta x \left(\sum_{j=1}^m \frac{u_{j-1}^n+u_j^n}2 +\sum_{j=m+1}^N \frac{v_{j-1}^n+v_j^n}2\right)\\
&=&\Delta x \left(\frac 12 u_0^n + \sum_{j=1}^{m-1} u_j^n +\frac 12 (u_m^n +v_m^n) +\sum_{j=m+1}^{N-1} v_j^n
+\frac 12v_N^n\right).
\end{eqnarray*}
This is the same formula as when $\overline{C}_n$ is obtained using \eqref{dnc2}. The
natural boundary conditions and mass lumping for the time derivative maintain the conservation property.

\subsubsection*{Discrete mass conservation for finite volume type coupling schemes}

In the case of the finite volume type scheme, we can simplify \eqref{dnc2} to
\begin{equation}
    \label{fv_cons}
C_n=u^n_1+\cdots +u^n_m+v^n_{m+1}+\cdots+v^n_{N}. 
\end{equation}
We see that the factors of $2$ are not needed,
so we obtain the conservation property using the discrete boundary conditions \eqref{onesided_bdry}
and heat flux coupling conditions \eqref{300} and \eqref{3001}. Again, we can substitute
any of the other fluxes for $J_{heat}(u_m^n,v_m^n)$.

The conservation property also extends to multidimensional equations. For finite volume schemes, this is straightforward.
The total computed quantity is just summed up of all cells. 
For finite elements on regular rectangular or quadrilateral meshes,
it is also clear how to define cells around the nodes. 
For Delauney triangulations, one can use the Voronoi cells around the nodes.
These cells are used in some finite volume methods, see e.g.\ Mishev \cite{MIS} and literature cited there.

\section{Numerical tests for the coupling conditions and the discrete conservation property}
\label{sec:tests}

In order to test numerical schemes, it is useful to have a problem with exact solutions. 
For the diffusion equation $w_t-Dw_{xx} =0$, one can easily generate solutions
by separation of variables, see. e.g.\ John \cite[Section 7.1]{bJOH}.
For us, an interesting family of solutions on the interval $[0,1]$ are
$$
w_n(x,t)=\e^{-D(n\pi)^2 t}\cos \left(n\pi x\right)+1
$$
for $n=1,2,3,\ldots$. They satisfy the homogeneous Neumann boundary
condition. We take $n=1$ then the initial data are $w(x)=\cos \left(\pi x\right)+1$ 
with initial total mass concentration $\int_{0}^{1} [\cos(\pi x)+1] \,dx=1$.
Due to the homogeneous Neumann boundary conditions, the total mass concentration 
remains unchanged for the exact solution. 
The solution will approach the constant $1$ asymptotically in time. So, 
total mass concentration is a good indicator for finding flaws in the code.
Using this solution, we validated our discretizations and the numerical Dirichlet-Neumann coupling conditions. 
The latter should reproduce the same exact solution up to rounding errors in a bi-domain computation. 
In our computations, this was indeed the case.

All our computations were done with {\sc Matlab} R2023b on the interval $[0,1]$ with double precision. 
We chose the spatial mesh size and set the time step to be in the relation
$\Delta t =\frac{0.4(\Delta x)^2}{\max (D_-,D_+)}$, i.e.\ $\nu_\pm\le 0.4$, 
to have a safety margin towards the stability bound at $0.5$.

\subsection{The coupling conditions}

\begin{figure}
\begin{subfigure}{.49\textwidth}
  \centering
  \includegraphics[width=1 \linewidth]{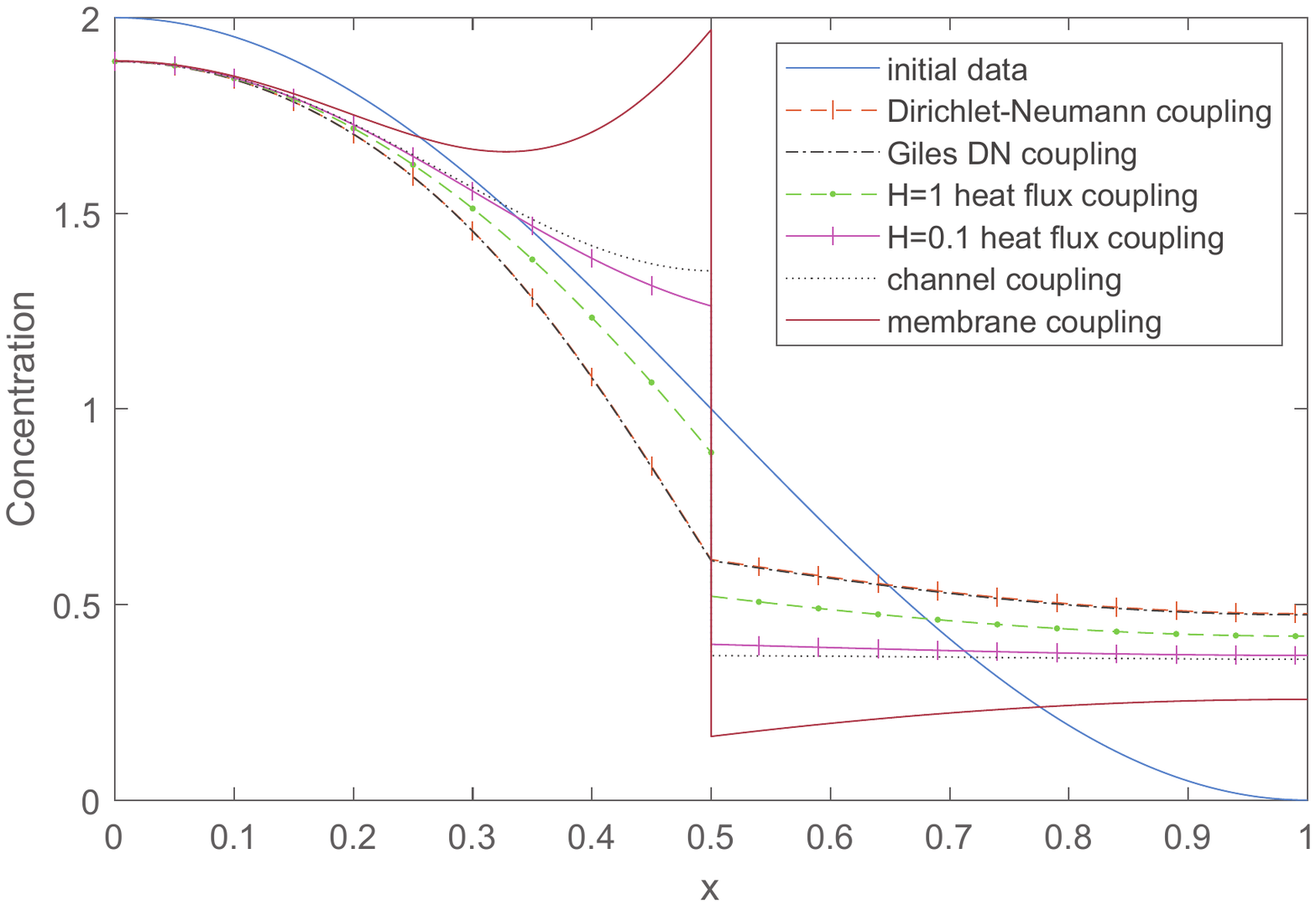}
  \caption{}
\end{subfigure}
\begin{subfigure}{.49\textwidth}
  \centering
  \includegraphics[width=1 \linewidth]{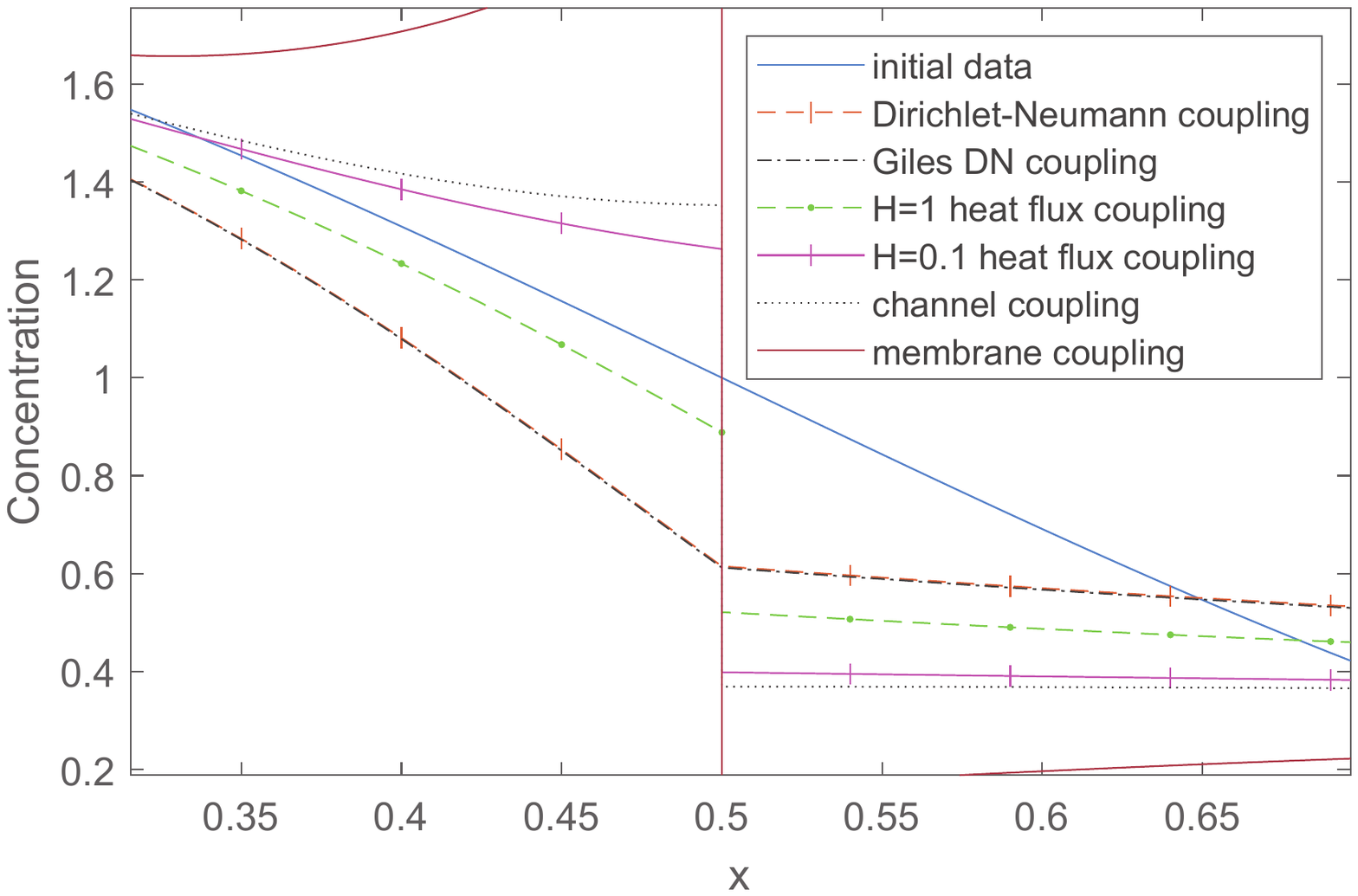}
  \caption{}
\end{subfigure}
\caption{\textbf{Bi-domain computations with six couplings:} 
\label{comp1}
The interval [0,1] is divided into 100 sub-intervals of length $\Delta x=0.01$. 
There were $3000$ time steps of length $\Delta t=4\cdot 10^{-5}$, giving a final time $T=0.12$. 
The figure on the right is a zoom into the coupling region.}
\end{figure}

The Dirichlet-Neumann coupling and the heat flux coupling are well established. 
The biophysical channel and pumping conditions are slightly
different. So, for comparison, we show how they act using a few examples.
We took the initial data $\cos(\pi x)+1$ as well as the diffusion coefficients $D_-=0.1$ and $D_+=1$.
Shown in Figure \ref{comp1} are the initial data and
six computations with the nodal based scheme. 

We used the coupling conditions \eqref{f16} and \eqref{giless1} together with 
the Dirichlet condition $u_m^{n+1}=v_m^{n+1}$ and then \eqref{2001}
with the heat flux coefficients $H=1$ and $H=0.1$. Further, we took the channel pumping with $\psi = 9.3954\text{E}-7$,
$\alpha = 1.497$, $\beta = 1.1949\text{E}-04$, $\gamma = 1.1556\text{E}-07$ and $\delta = 1.1444\text{E}-07$. 
For the membrane pumping we chose $P_l=0.02$, $P_p=1$ and $K_d=0.2$.
The values were almost all taken from a table in Thul and Falcke \cite{l11}. 
These authors require $\psi = 9.3954$ and $P_p =40$. The values correspond to a larger heat flux.
But this leads to a bad scaling within the figure. So, we modified them for the sake of our comparison.
The higher values were used for local pumping in a 3 dimensional model where the concentration spreads out and becomes more diluted.
A one dimensional model needs only a smaller value for the coefficients in order to produce the same effect.

We made the spatial mesh rather course with $\Delta x=0.01$ in order to highlight differences in the couplings.
In the zoom, we cut off the extreme values of the membrane coupling to show the other
solutions more clearly. To the left of $x=0.5$ from top to bottom, the solutions are with
membrane, then channel coupling, the initial data, heat flux coupling with $H=0.1$, then $H=1$, and finally Dirichlet-Neumann and Giles coupling.
The latter two are continuous. To the right, the order of the discontinuous solutions is reversed
and they all are below the continuous couplings.

We see in the computational results that the solution approaches a quasi constant state 
faster on the right sub-interval due to the larger diffusion coefficient. 
Zooming in even further into the Graph of the Dirichlet-Neumann coupling would show
that the Giles formula does deviate from the correct coupling.
The nodal values at $x=0.5$ differ by $2.730344211099967\text{E}-03$.
The smaller heat flux coefficient of $H=0.05$ models a less permeable membrane than the case $H=1$. 
So, the overall equilibrium value $u(x)=v(x)=1$ is approached much slower. 

The pumping goes from right to left, increasing the concentration on the left hand side 
at the expense of the right hand side. 
The small value of $H=0.1$ represents a slower flux through the membrane
at $x=0.5$ than with $H=1$.
So, the left-hand side is loosing mass at a slower rate.

\begin{figure}
\begin{subfigure}{.49\textwidth}
  \centering
  \includegraphics[width=1 \linewidth]{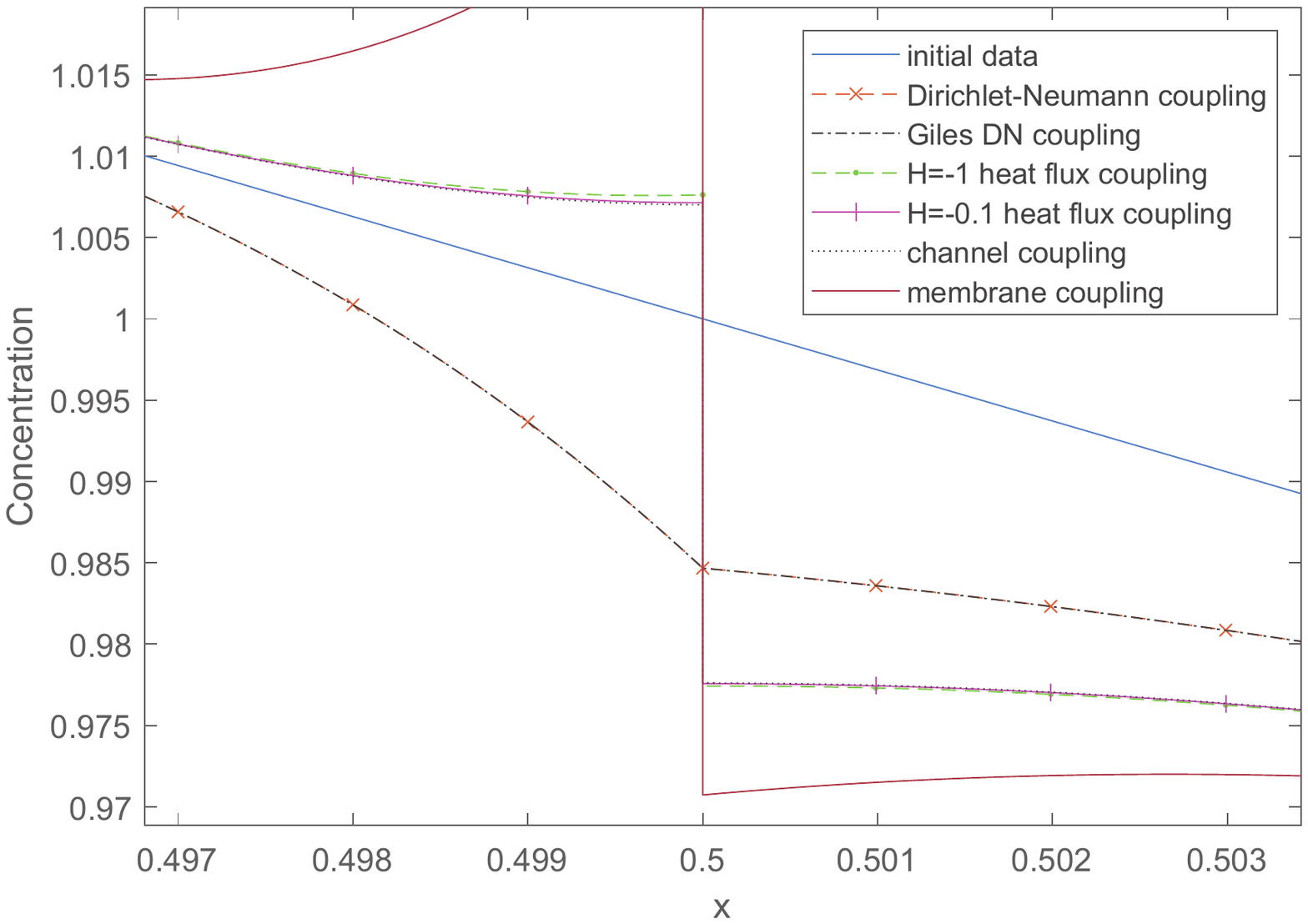}
  \caption{}
\end{subfigure}
\begin{subfigure}{.49\textwidth}
  \centering
  \includegraphics[width=1 \linewidth]{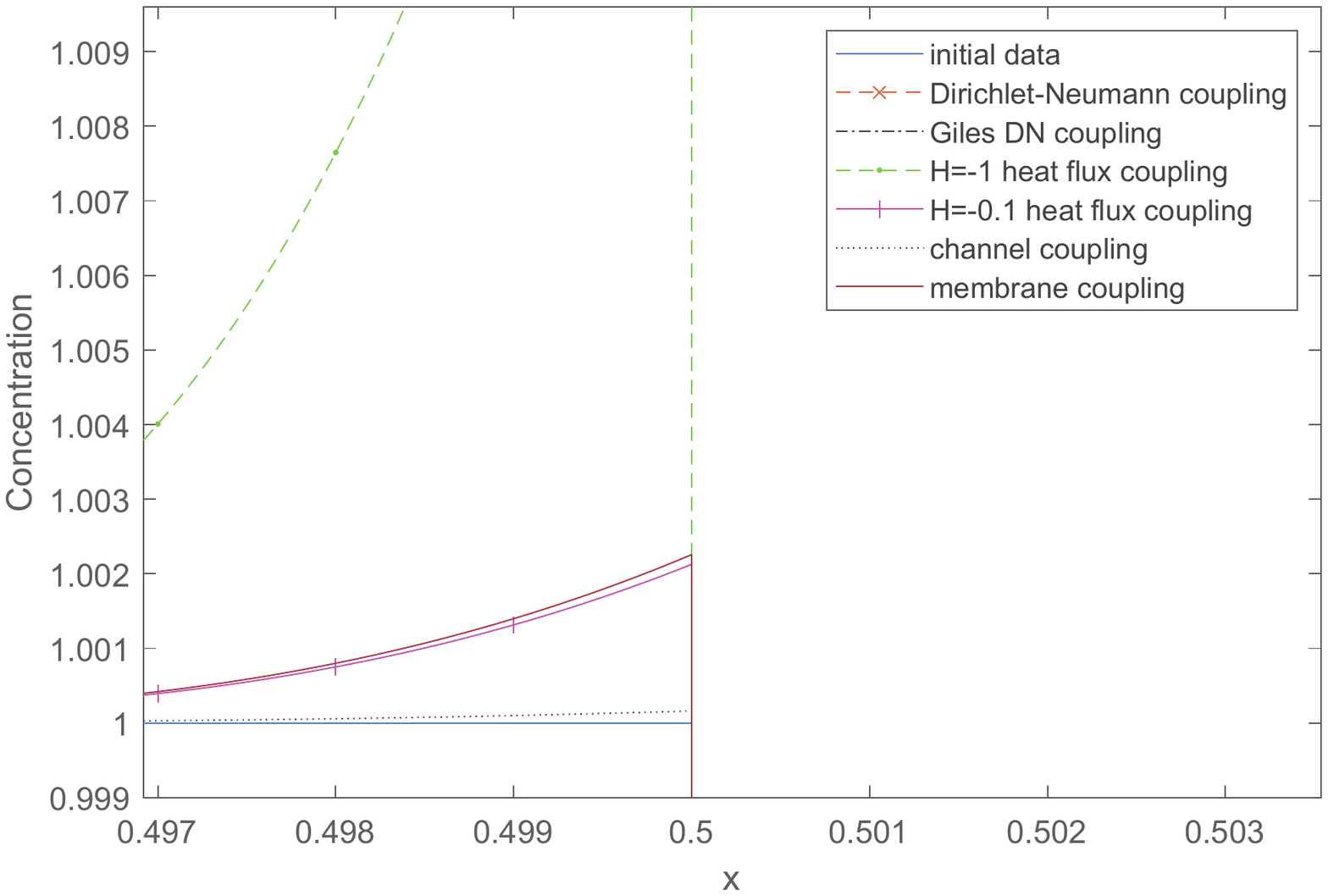}
  \caption{}
\end{subfigure}
\caption{\textbf{Bi-domain computations with negative heat fluxes:} 
\label{comp5}
The interval [0,1] is divided into 100,000 sub-intervals of length $\Delta x=0.00001$. 
There were $10^6$ time steps of length $\Delta t=4\cdot 10^{-11}$, giving a final time $T=0.00004$. The figure on the left 
is a zoom into the coupling region for the cosine initial data.
The figure on the right is a zoom into the upper part of the coupling region
for the piecewise constant initial data.}
\end{figure}

For comparison of the influence of the heat flux coefficients in the coupling, 
we also made two computations with negative heat flux coefficients using $H=-1$, $H=-0.01$, $\psi = -9.3954\text{E}-7$
and $P_l=-0.02$. Figure \ref{comp5} shows two zooms into the solutions.
Negative heat flux coefficients mean that there is no natural heat or concentration flux 
against the gradient but an active flux mechanism working
in the direction of the gradient. We used the finer mesh, which we will be using below. The first computation 
is with the above cosine initial data. The second has piecewise constant initial data
with $u(x)=1$ and $v(x)=0.06$, cp.\ Figure \ref{comp3}. In Figure \ref{comp5}(a), to the left of the interface from top
to bottom are the membrane coupling, $H=-1$, $H=-0.1$, and channel coupling solutions. The latter three are quite
close to each other. In Figure \ref{comp5}(a), we have the order of solutions $H=-1$, membrane coupling, $H=-0.1$, and channel
coupling. The latter is close to the initial data. We clearly see that in each case more heat or mass 
is moved from right to left as compared to the other computations with positive coefficients.

\begin{figure}
\begin{subfigure}{.49\textwidth}
  \centering
  \includegraphics[width=1 \linewidth]{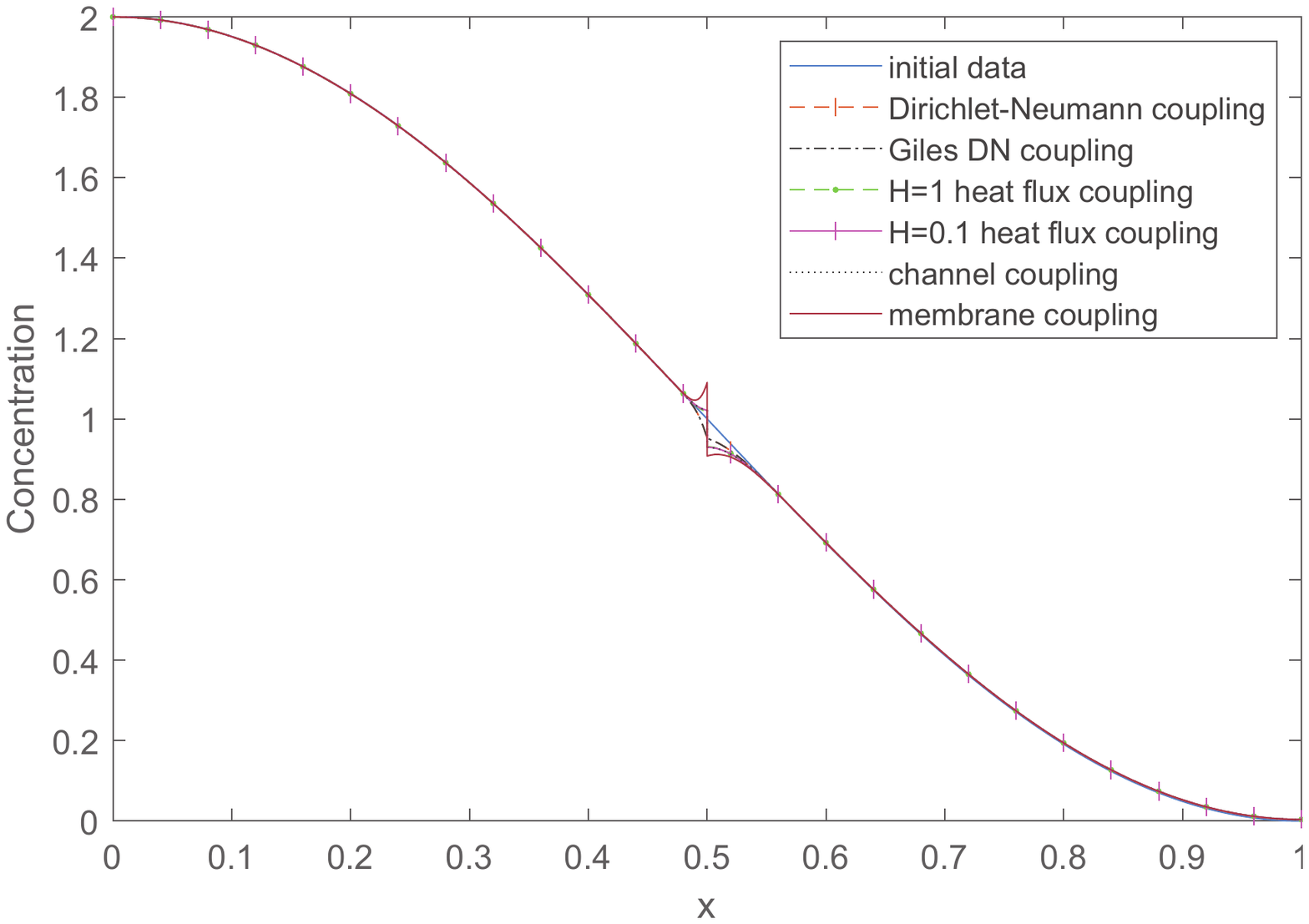}
  \caption{}
\end{subfigure}
\begin{subfigure}{.49\textwidth}
  \centering
  \includegraphics[width=1 \linewidth]{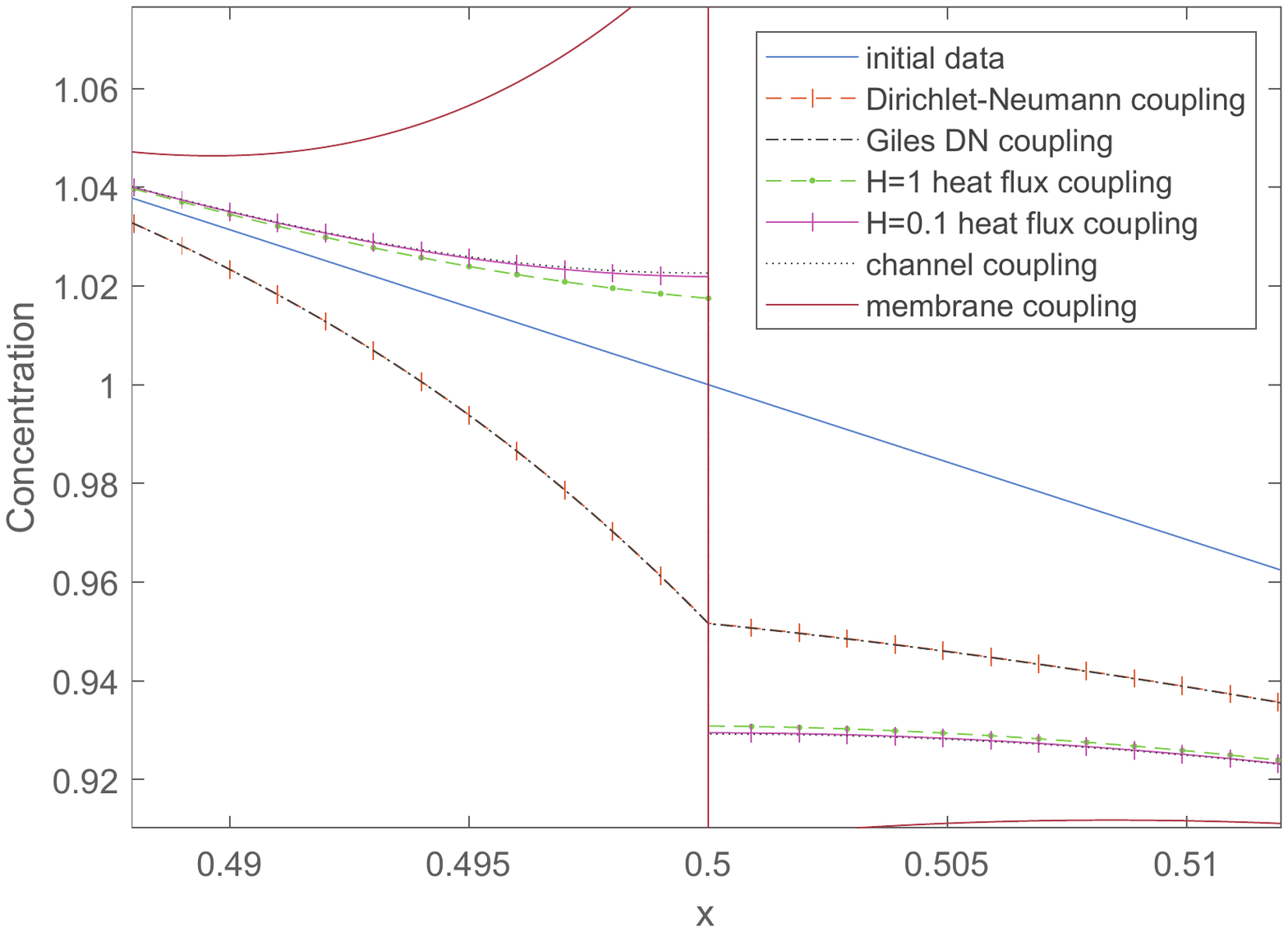}
  \caption{}
\end{subfigure}
\caption{\textbf{Bi-domain computations with six couplings:} 
\label{comp2}
The interval [0,1] is divided into 10,000 sub-intervals of length $\Delta x=0.0001$. 
There were $10^5$ time steps of length $\Delta t=4\cdot 10^{-9}$, giving a final time $T=0.0004$. 
The figure on the right is a zoom into the coupling region.}
\end{figure}

\subsection{Discrete mass conservation of coupling conditions}

We made the same computation as in the previous section on a finer mesh with $\Delta x = 0.0001$ and $10^5$ time steps 
of size $\Delta t=4\cdot 10^{-9}$ in order to
demonstrate the mass conservation property. 
The solutions are shown in Figure \ref{comp2}. 
The extreme values at $x=0.5$ are taken by the membrane pumping solution. 
The channel pumping solution is on the left hand side only slightly above and on the right hand side
below the $H=0.1$ heat flux solution. Then come the $H=1$ heat flux solution 
and the continuous values of the initial data.
Due to the resolution of the figure, the graph of the
Dirichlet-Neumann coupling solutions are hidden behind the one using 
the Giles coupling. They are optically quite close. 
But, their difference in solution value at $x=0.5$ is $8.133493461626173\text{E}-05$.

We had a discrete initial total mass concentration of
$C(0)=1.000000000000001$. Then, we obtained the final values 
$C(T)=0.9999975790344722$ for the correct Dirichlet-Neumann coupling and 
$C(T)=0.9980606224422512$ for the Giles coupling. All of these were computed using \eqref{dnc2}.
The computational error of the correct coupling was $|C(T)-C((0)|=4.432343381211012\text{E}-12$. 
We clearly see that the Giles coupling produces 
a much larger error. It is $|C(T)-C(0)|=2.420965528937558\text{E}-06$ by losing mass.
The heat flux couplings gave $C(T) = 0.9999999999955704$
and $C(T)= 9.999999999955695$. For the pumping the values were $C(T) = 0.9999999999955702$ 
and $C(T) = 0.9999999999955707$. Except for the Giles coupling, the values 
are very reasonable deviations of order $10^{-12}$ from the exact value $1$.

\begin{figure}
\begin{subfigure}{.49\textwidth}
  \centering
  \includegraphics[width=1 \linewidth]{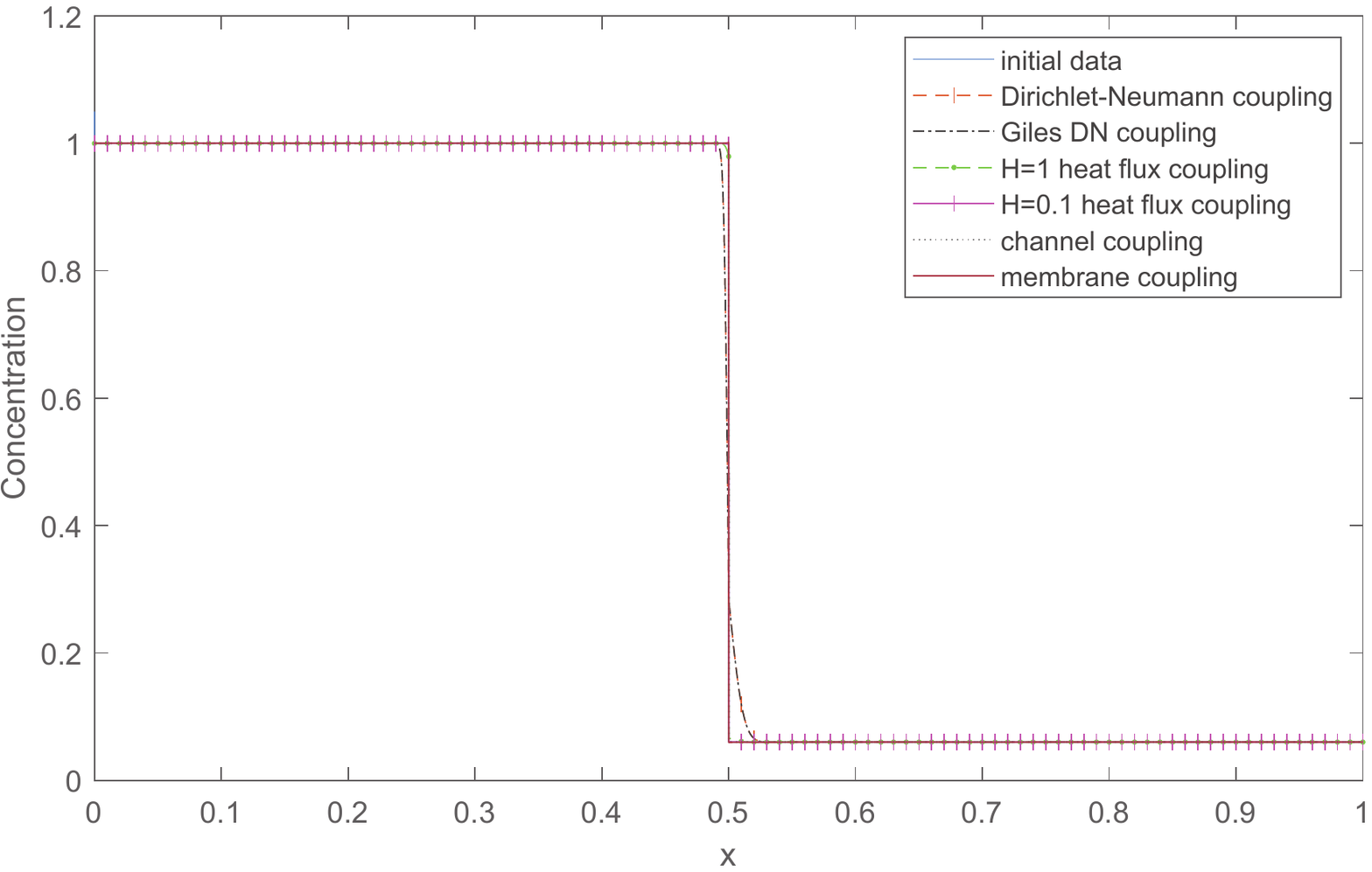}
  \caption{}
\end{subfigure}
\begin{subfigure}{.49\textwidth}
  \centering
  \includegraphics[width=1 \linewidth]{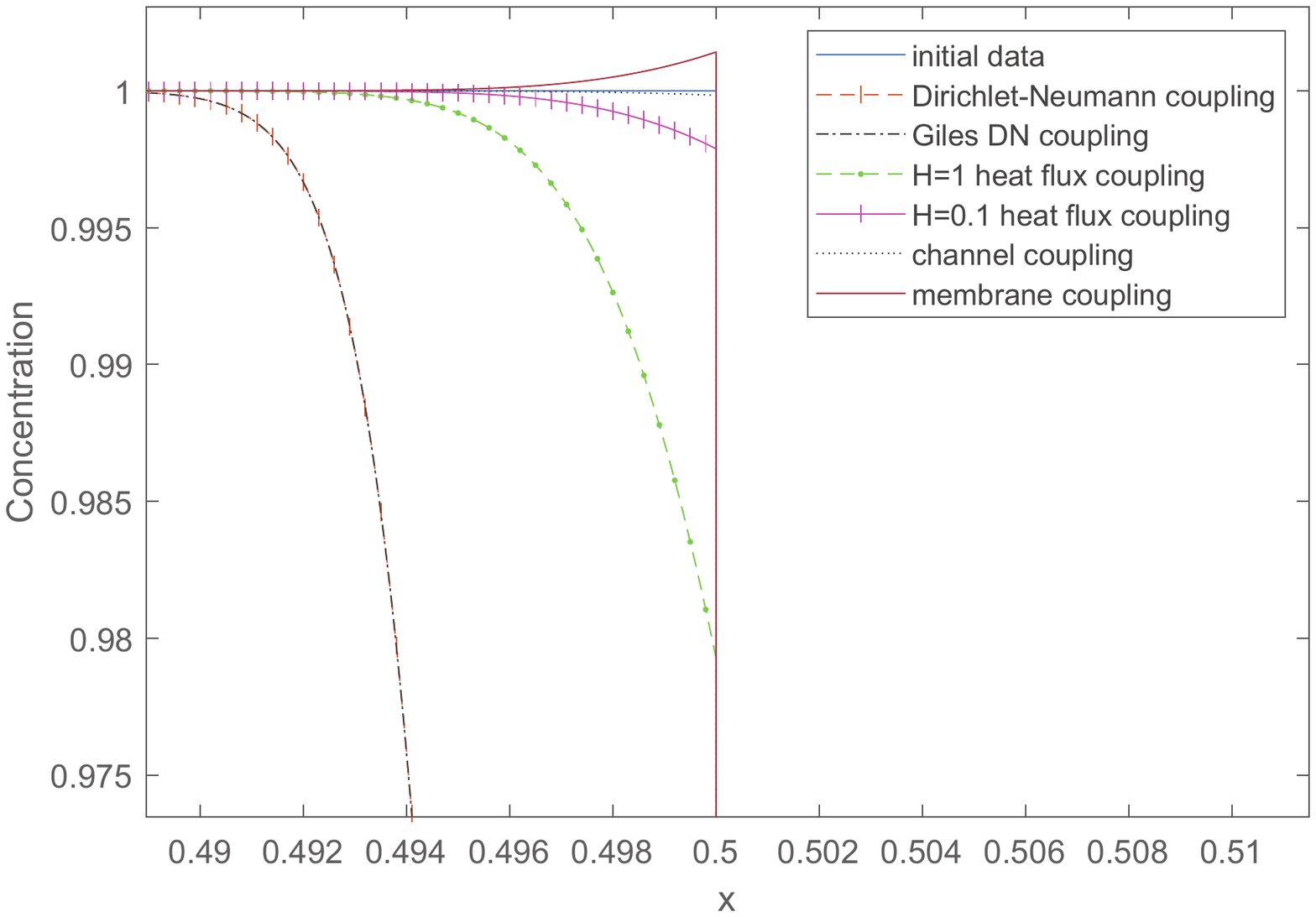}
  \caption{}
\end{subfigure}
\caption{\textbf{Bi-domain computations with six couplings:} 
\label{comp3}
The interval [0,1] is divided into 100,000 sub-intervals of length $\Delta x=0.00001$. 
There were $10^6$ time steps of length $\Delta t=4\cdot 10^{-11}$, giving a final time $T=0.00004$. 
The figure on the right is a zoom into the upper part of the coupling region.}
\end{figure}

Next, we took piecewise constant initial data with $u(x)=1$ and $v(x)=0.06$, see Figure \ref{comp3}.
The exact initial total concentration is $0.53$. The initial discretized total 
concentration was $C(0)=0.5300000000000005$. 
For the Dirichlet-Neumann coupling we obtained $C(T)=0.5300009399993384$
and for the Giles coupling $C(T)=0.5299973682972777$. The results for 
the heat fluxes were $C(T)=0.5299999999993753$, $C(T)=0.5299999999994137$
and for the pumping couplings $C(T)=0.5299999999994603$, $C(T)=0.5299999999993950$.
The Dirichlet-Neumann coupling has a larger computational error of $|C(T)-C(0)|=9.399993379233251\text{E}-07$
then the other correct couplings, which are of the order $10^{-13}$. 
This seems to be caused initially by the Dirichlet condition
which does not do well with the initial discontinuity. With an error
of $|C(T)-C(0)|=2.631702722744045\text{E}-06$, the error in the Giles coupling is not
so pronounced in the example.

We did the same computation taking the non-conservative one-sided differences
giving \eqref{300} and \eqref{3001} for the flux couplings. We obtained for the heat fluxes
$C(T)=0.5300000706650079$, $C(T)=0.5300000706650079$ and for the pumping fluxes
$C(T)=0.5300000005493521$, $C(T)=0.5299999951641959$. The computational error
$|C(T)-C(0)|$ is considerably higher with the non-conservative differences. The order
is $10^{-8}$ for the heat fluxes, $10^{-10}$ for the channel pumping flux and
$10^{-9}$ for the membrane pumping flux. The plots are not much different from Figure \ref{comp3}.

\subsection*{The finite volume type discretization}
\begin{figure}
\begin{subfigure}{.49\textwidth}
  \centering
  \includegraphics[width=1 \linewidth]{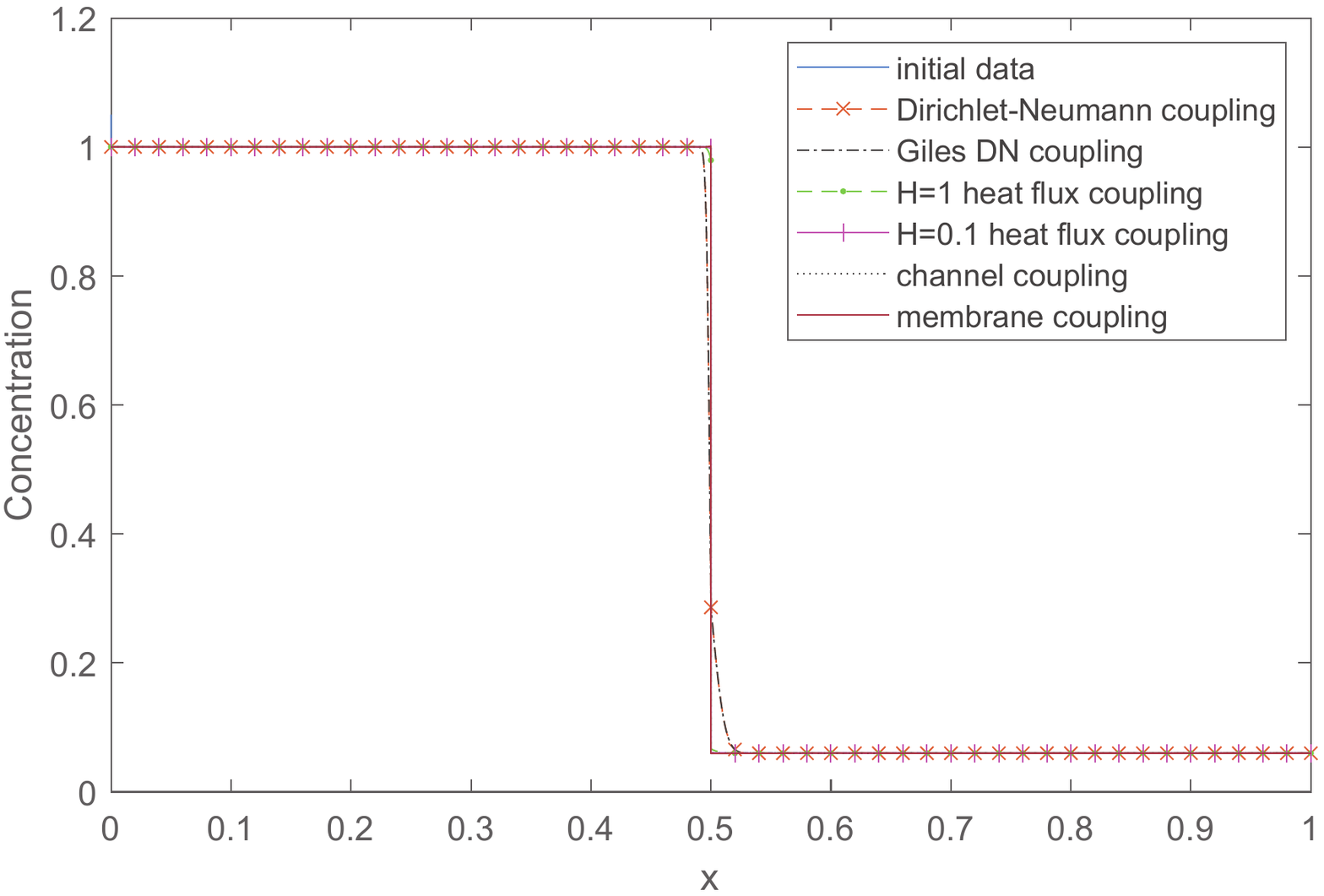}
  \caption{}
\end{subfigure}
\begin{subfigure}{.49\textwidth}
  \centering
  \includegraphics[width=1 \linewidth]{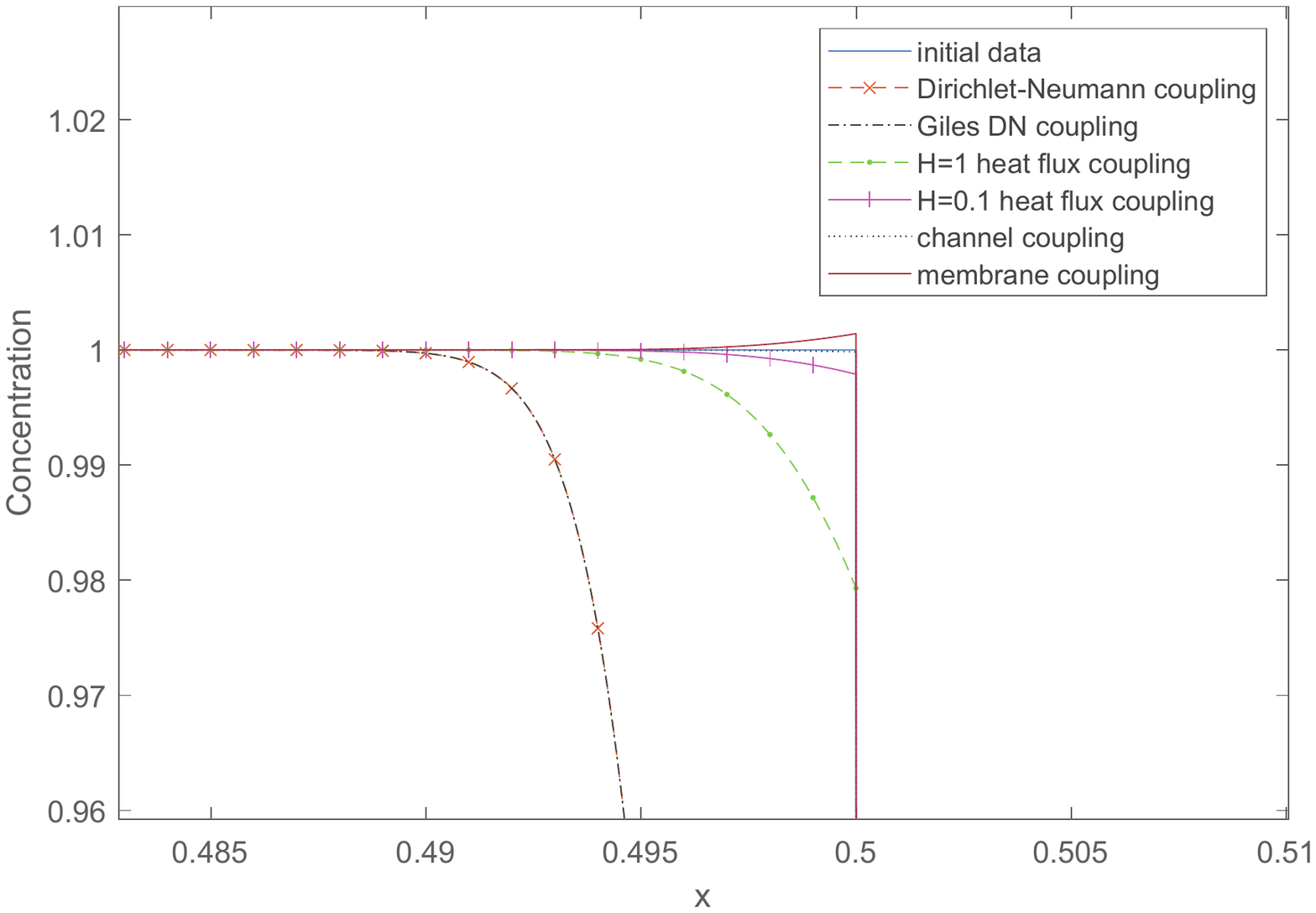}
  \caption{}
\end{subfigure}
\caption{\textbf{Finite volume type bi-domain computations with six couplings:} 
\label{comp4}
The interval [0,1] is divided into 100,000 sub-intervals of length $\Delta x=0.00001$. 
There were $10^6$ time steps of length $\Delta t=4\cdot 10^{-11}$, giving a final time $T=0.00004$. 
The figure on the right is a zoom into the upper part of the coupling region.}
\end{figure}

Finally, we redid the preceding computations with the piecewise constant initial data 
using the same parameters with the finite volume type discretization
that was used for the nodal based scheme, see Figure \ref{comp4}.
The nodal points were shifted, the mass conserving flux coupling formulas
\eqref{300} and \eqref{3001} as well as boundary condition
formulas \eqref{onesided_bdry} were used. Further, the discrete total mass concentration was
computed using \eqref{fv_cons}. The initial discrete total mass concentration was again $C(0)=0.5300000000000005$.
The Dirichlet-Neumann coupling gave a much better result of $C(T)=0.5299999999993383$ with an error of
$|C(T)-C(0)|=6.621370118864434\text{E}-13$. Here, we see that the node at the initial discontinuity caused the large
error in the nodal based scheme. The Giles coupling coupling gives 
$C(T)= 0.5299964295714319$ with an error of  $|C(T)-C(0)|=3.570428568577810e-06$. 
Again, it is much larger than the correct coupling. The heat flux couplings give $C(T)=0.5299999999993757$,
$C(T)=0.5299999999994144$ and the pumping couplings $C(T)=0.5299999999994608$,
$C(T)=0.5299999999993955$. All computational errors are of size $10^{-13}$
as before.

Some further computational results for the conservation property were given in Munir \cite[Section 5.5]{MUN}. 

\subsection{The homogeneous Neumann boundary condition}

For the discrete mass conservation, let us first look at the homogeneous Neumann boundary conditions
\eqref{onesided_bdry} and \eqref{central_bdry}. We mentioned that the
former is good with the finite volume type scheme and the latter with the nodal-based scheme. 
The error when taking the non-conservative discretization at
the boundary is most pronounced when the solution has a steep gradient
at the boundary. We want to show this effect.

We did some computations for a simple single-domain problem with the nodal-based scheme with $D=0.001$. 
As initial data we took the function $u(x) = 100\sqrt{x(1-x)}$ for $x\in [0,1]$
with the initial total mass concentration $C=100\int_{0}^{1} \sqrt{x(1-x)}\,dx=39.269908169872415$.

First, we took the boundary condition \eqref{central_bdry}.
We ran three computations, the first with $\Delta x = 0.01$ and 500 time steps
of size $\Delta t = 0.04$ to the final time $T=20$.
The total approximated initial concentration was $C(0)=39.22835638873124$
with a large initial dicretization error of $|C(T)-C(0)|= 4.155178114117319\text{E}-02$. As the final
total concentration, we obtained $C(T) =39.22835638873113$. Only the last two decimals
are changed, giving a computational error of $|C(T)-C(0)|=1.136868377216160\text{E}-13$. Here, the initial error heavily outweighs the computational error.

Next we used $\Delta x = 0.001$ with $50.000$ time steps of size $\Delta t = 0.0004$ to
the final time $T=20$.
We obtained $C(0) = 39.26859346252676$ with an initial discretization error of $|C(T)-C(0)|=1.314707345656529\text{E}-03$. The final total concentration was
$C(T)=39.26859346251628$. 
\text{E}ven though the initial total concentration is
more precise due to the finer spatial mesh, we have a higher computational error of
$|C(T)-C(0)|=1.048050535246148\text{E}-11$. We see the effect of error accumulation due to the
many time steps. 

We also used a much finer mesh of $\Delta x = 10^{-7}$ to get 
a better discretization of the initial values.
We had a time step of $\Delta t = 4\cdot 10^{-12}$ due to the stiffness of the problem and took $10^5$ time steps to the
final time $T=4\cdot 10^{-7}$. We have the initial total
concentration $C(0)=39.26990816855763$ with an initial discretization error of $|C-C(0)|=1.314788278250489e-09$ and the final
total concentration $C(T)=3.926990816855260$ with a computational error of $|C(T)-C(0)|=5.023537141823908\text{E}-12$.

We repeated these computations with the one-sided differences \eqref{onesided_bdry} at both boundaries. The initial errors were the same.
In the first case, we have $C(T)=38.91134647144038$, giving a computational error of $|C(T)-C(0)|=0.3170099172908607$, 
which is very large compared to the result with the conservative boundary condition.
On the finer mesh, we obtained $C(T)=39.23609630251410$ with a slightly
smaller computational error of $|C(T)-C(0)|= 0.03249716001266023$ due to the finer mesh. On the very fine mesh 
at the smaller final time, we obtained $C(T)=39.26990812490996$ with
a computational error $|C(T)-C(0)|=4.364766681419496\text{E}-08$. It is clearly larger than in the case of the conservative boundary condition.

\section{Conclusion}

We introduced an improved concept of the discrete conservation property for finite difference schemes.
It was demonstrated how this property is important to determine the correct choice of
flux boundary and coupling conditions for the diffusion or heat equation. This was done both
theoretically and using a number of computational examples. The violation of this property
can produce significant computational errors when there are steep gradients near a boundary
or a coupling interface. Our study examined various coupling conditions, including well-established methods such as Dirichlet-Neumann coupling, heat flux coupling, and specific channel and pumping flux conditions from biophysics. Additionally, we highlighted the crucial differences in handling fluxes between nodal-based and finite volume schemes, providing clarity through multiple examples discussed in this work.

\subsection*{Acknowledgements} The first author would like to thank the government of Pakistan and the
German Academic Exchange Service (DAAD) for generously supporting his doctoral research. The authors are also
grateful to Martin Falcke for cooperating with them on the topic of modeling and computing 
calcium dynamics in cells with reaction-diffusion systems.
This joint work provided the impetus for this paper.

\subsection*{Data Availability}
Data sharing is not applicable to this article as no datasets were generated or analyzed during the current study.

\subsection*{ Conflicts of Interest}
The authors declare no conflicts of interest.

\bibliography{book.bib,tajrefrence.bib}

\begin{thebibliography}{10}

\bibitem{l116}
{\sc E.~Carr and N.~March}, {\em Semi-analytical solution of multilayer
  diffusion problems with time-varying boundary conditions and general
  interface conditions}, Applied Mathematics and Computation, 333 (2018),
  pp.~286--303.

\bibitem{l9}
{\sc N.~Chamakuri}, {\em Adaptive numerical simulation of reaction-diffusion
  systems}, PhD thesis, University of Magdeburg, 2007.

\bibitem{errera2013optimal}
{\sc M.-P. Errera and S.~Chemin}, {\em Optimal solutions of numerical interface
  conditions in fluid--structure thermal analysis}, Journal of Computational
  Physics, 245 (2013), pp.~431--455.

\bibitem{l8}
{\sc M.~Falcke}, {\em On the role of stochastic channel behavior in
  intracellular {Ca2+} dynamics}, Biophysical journal, 84 (2003), pp.~42--56.

\bibitem{GIL}
{\sc M.~B. Giles}, {\em Stability analysis of numerical interface conditions in
  fluid-structure thermal analysis}, Internat. J. Numer. Methods Fluids, 25
  (1997), pp.~421--436.

\bibitem{chand}
{\sc W.~D. Henshaw and K.~K. Chand}, {\em A composite grid solver for conjugate
  heat transfer in fluid--structure systems}, Journal of Computational Physics,
  228 (2009), pp.~3708--3741.

\bibitem{bHUVE}
{\sc W.~Hundsdorfer and J.~Verwer}, {\em Numerical Solution of Time-Dependent
  Advection-Diffusion-Reaction Equations}, vol.~33 of Springer Series in
  Computational Mathematics, Springer-Verlag, Berlin, 2003.

\bibitem{b8}
{\sc J.~C. Jaeger and H.~S. Carslaw}, {\em Conduction of Heat in Solids},
  Clarendon Press Oxford, 1959.

\bibitem{bJOH}
{\sc F.~John}, {\em Partial Differential Equations}, Springer-Verlag, New
  York-Heidelberg-Berlin, 1978.

\bibitem{bJOHN}
{\sc C.~Johnson}, {\em Numerical Solution of Partial Differential Equations by
  the Finite Element Method}, Cambridge University Press, Cambridge - New York,
  1987.

\bibitem{lemarie2015analysis}
{\sc F.~Lemari{\'e}, E.~Blayo, and L.~Debreu}, {\em Analysis of
  ocean-atmosphere coupling algorithms: consistency and stability}, Procedia
  Computer Science, 51 (2015), pp.~2066--2075.

\bibitem{MIS}
{\sc I.~D. Mishev}, {\em Finite volume methods on {V}oronoi meshes}, Numer.
  Methods Partial Differential Equations, 14 (1998), pp.~193--212.

\bibitem{moretti2018stability}
{\sc R.~Moretti, M.-P. Errera, V.~Couaillier, and F.~Feyel}, {\em Stability,
  convergence and optimization of interface treatments in weak and strong
  thermal fluid-structure interaction}, International Journal of Thermal
  Sciences, 126 (2018), pp.~23--37.

\bibitem{bMOMA}
{\sc K.~W. Morton and D.~F. Mayers}, {\em Numerical Solution of Partial
  Differential Equations}, Cambridge University Press, Cambridge, second~ed.,
  2005.
\newblock An introduction.

\bibitem{MUN}
{\sc T.~Munir}, {\em Analysis of coupling interface problems for bi-domain
  diffusion equations}, PhD thesis, Otto von Guericke University, Magdeburg,
  2020.

\bibitem{CMW2}
{\sc T.~Munir, N.~Chamakuri, and G.~Warnecke}, {\em On conservative, stable
  boundary and coupling conditions for diffusion equations {I}{I} -
  {S}tability}.
\newblock Manuscript in preparation.

\bibitem{CMW3}
\leavevmode\vrule height 2pt depth -1.6pt width 23pt, {\em On conservative,
  stable boundary and coupling conditions for diffusion equations {I}{I}{I} -
  {T}he conservation property for varying mesh sizes and implicit schemes}.
\newblock Manuscript in preparation.

\bibitem{bQUA}
{\sc A.~Quarteroni}, {\em Numerical Models for Differential Problems}, vol.~16
  of MS\&A. Modeling, Simulation and Applications, Springer, Cham, 2017.

\bibitem{bQUVA1}
{\sc A.~Quarteroni and A.~Valli}, {\em Numerical Approximation of Partial
  Differential Equations}, vol.~23 of Springer Series in Computational
  Mathematics, Springer-Verlag, Berlin, 1994.

\bibitem{bQUVA}
{\sc A.~{Quarteroni and A. Valli}}, {\em Domain Decomposition Methods for
  Partial Differential Equations}, Numerical Mathematics and Scientific
  Computation, The Clarendon Press, Oxford University Press, New York, 1999.
\newblock Oxford Science Publications.

\bibitem{roe2007stability}
{\sc B.~Roe, A.~Haselbacher, and P.~H. Geubelle}, {\em Stability of
  fluid--structure thermal simulations on moving grids}, International journal
  for numerical methods in fluids, 54 (2007), pp.~1097--1117.

\bibitem{roe2008combined}
{\sc B.~Roe, R.~Jaiman, A.~Haselbacher, and P.~H. Geubelle}, {\em Combined
  interface boundary condition method for coupled thermal simulations},
  International journal for numerical methods in fluids, 57 (2008),
  pp.~329--354.

\bibitem{bTHO}
{\sc J.~Thomas}, {\em Numerical Partial Differential Equations: Finite
  Difference Methods}, vol.~22 of Texts in Applied Mathematics,
  Springer-Verlag, New York, 1995.

\bibitem{bTHOME}
{\sc V.~Thom\'{e}e}, {\em Galerkin Finite Element Methods for Parabolic
  Problems}, vol.~25 of Springer Series in Computational Mathematics,
  Springer-Verlag, Berlin, second~ed., 2006.

\bibitem{l10}
{\sc R.~Thul}, {\em Analysis of intracellular reaction diffusion systems: The
  stochastic medium Calcium}, PhD thesis, Free University, Berlin, November,
  2004.

\bibitem{l11}
{\sc R.~Thul and M.~Falcke}, {\em Release currents of {IP3} receptor channel
  clusters and concentration profiles}, Biophysical journal, 86 (2004),
  pp.~2660--2673.

\bibitem{bTOWI}
{\sc A.~{Toselli and O. Widlund}}, {\em Domain Decomposition
  Methods---Algorithms and Theory}, vol.~34 of Springer Series in Computational
  Mathematics, Springer-Verlag, Berlin, 2005.

\bibitem{bWAMI}
{\sc R.~Wait and A.~R. Mitchell}, {\em Finite Element Analysis and
  Applications}, A Wiley-Interscience Publication, John Wiley \& Sons, Inc.,
  New York, 1985.

\bibitem{zhang}
{\sc H.~Zhang, Z.~Liu, E.~Constantinescu, and R.~Jacob}, {\em Stability
  analysis of interface conditions for ocean-atmosphere coupling}, J. Sci.
  Comput., 84 (2020), pp.~Paper No. 44, 25.

\end{thebibliography}
\bibliographystyle{siam}

\end{document}